%% file: MT.tex
\documentclass[11pt, a4paper, UKenglish, amsthm]{article}

\usepackage[UKenglish]{babel}
\usepackage{MA_Titlepage}

\usepackage{graphicx}
\usepackage{amsfonts, ifthen, amsmath, amssymb, esint}
\usepackage{ntheorem}

\usepackage{skak}
\usepackage{enumitem}
\usepackage[style=alphabetic, backend=bibtex8, doi=false, isbn=false]{biblatex}
\bibliography{LiteratureMT}

\newcommand{\grad}{\ensuremath{\nabla}}
\newcommand{\ve}{\ensuremath{\varepsilon}}
\newcommand{\vp}{\ensuremath{\varphi}}
\newcommand{\al}{\ensuremath{\alpha}}
\newcommand{\ga}{\ensuremath{\gamma}}
\newcommand{\ta}{\ensuremath{\tau}}
\newcommand{\ot}{\ensuremath{\otimes}}
\newcommand{\xh}{\ensuremath{\hat{x}}}
\newcommand{\xp}{\ensuremath{\hat{x}^{\perp}}}
\newcommand{\xa}{\ensuremath{|x|}}
\newcommand{\Id}{\ensuremath{Id_{2\times2}}}
\newcommand{\xpp}{\ensuremath{{\hat{x}}^{\perp}\ot\hat{x}}^{\perp}}
\newcommand{\xhh}{\ensuremath{{\hat{x}}\ot\hat{x}}}

\newcommand{\IvK}{\ensuremath{I_{h,\Delta}^{\mbox{{vK}}}(U,W)}}
\newcommand{\kvK}{\ensuremath{\kappa_W^{vK}}}
\newcommand{\wt}[1]{\widetilde{#1}}

\newcommand{\R}{\mathbb{R}}

\newcommand{\N}{\mathbb{N}}
\newcommand{\bb}[1]{\ensuremath{{\mathbb{#1}}}}
\newcommand{\bbs}[1]{\ensuremath{{\mathbb{R}^{#1}}}}
\newcommand{\A}{\mathcal{A}}
\newcommand{\B}{\mathcal{B}}
\newcommand{\La}{\mathcal{L}}

\newcommand{\WSob}{\ensuremath{W^{1,2}([0,1])}}
\newcommand{\Ite}{\ensuremath{I_{\ta,\ve}}}
\newcommand{\osc}[2]{\ensuremath{\underset{{#1}}{{\mbox{osc}}}\{#2\}}}
\newcommand{\ls}{\ensuremath{\lesssim}}
\newcommand{\gs}{\ensuremath{\gtrsim}}
\usepackage[capitalize]{cleveref}

\newtheorem{de}{Definition}[section]
\newtheorem{thm}[de]{Theorem}

\newtheorem{lem}[de]{Lemma}

\allowdisplaybreaks
\authornew{Marcel Dengler}
\geburtsdatum{16th February 1991}
\geburtsort{Balingen, Germany}
\date{\today}

\betreuer{Advisor: Dr. Heiner Olbermann}
\zweitgutachter{Second Advisor: Prof. Dr. Sergio Conti}
\institut{Institute for Applied Mathematics}
\title{Variational problems in thin elastic structures}
\ausarbeitungstyp{Master's Thesis  Mathematics}

\begin{document}

\maketitle

\tableofcontents
\include{MTIntro}

\include{MT11}

\include{MT11LB}


\include{MT21}

\include{MT22}

\printbibliography
\addcontentsline{toc}{section}{References}

\end{document}

%% file: MTIntro.tex
{\textbf{\Large{Acknowledgements}}}\\
\addcontentsline{toc}{section}{Acknowledgements}
I want to express my gratitude to Heiner Olbermann for introducing me to a wide range of the subject. Also I would like to thank him for the time and energy he invested and especially for the numerous Skype sessions which always resulted in new deep thoughts.\\
I'm thankful to Sergio Conti for cosupervising my thesis and inspiring discussions leading to new perspectives.  

\begin{flushright}
Marcel Dengler\\
D\"usseldorf, \date{\today}
\end{flushright}
\newpage

\section{Introduction}
In the last few years there has been a lot of interest in the area of crumpling and buckling thin elastic sheets. A key point to the understanding of such structures are the so called energy scaling laws, i.e. estimates on the energy which show it's behavior depending on the `thickness' $h$ of the sheet and other important parameters. The sheets being thin allows among other things the usage of 2-dimensional plate models, i.e. they can be seen as a two dimensional surface. Two dimensional models are approximations at the real world. The 3D-model is more complicated. For a relation between those models and a very useful tool to establish energy scaling laws can be found in the famous paper by James, Friesecke and Müller \autocite{MJF02}. Even if we only consider thin structures the thickness $h$ still enters as a parameter in the energy. We will think of $h$ as being very small.\\

When it comes to crumpling of sheets there are in particular two interesting situations, namely ridges and singularities. While this work deals with the latter one, a discussion of a minimal ridge in the Föppl-von Kármán model can be found in \autocite{V02}. 
A nice overview to the crumpling conjecture, together with a proof of the upper bound on the energy and results for special cases of the lower bound can be found in \autocite{CM08}. Note that Conti and Maggi prove their results for a 2D- and the general 3D-model.\\


The thesis consists of two independent parts the general theme being energy scaling laws of thin sheets.\\

The first part is related to \autocite{Ol15, Ol16, MO14}. In \autocite{Ol15} it is shown that there is an optimal energy scaling law for the geometry of a singular cone in different models. We are mainly interested in the Föppl-von Kármán model therefore we state the main result achieved in \autocite{Ol15} only in this picture, see \autocite[Thm 3]{Ol15}. The result is
\begin{equation}
2\pi\Delta^2h^2(\log{\frac{1}{h}}-2\log|\log h|-C)\le \min I_{h,\Delta}^{vK}\le 2\pi\Delta^2h^2(\log{\frac{1}{h}}+C).
\label{eq:OlB1}
\end{equation}
where $0<\Delta<1$ and $h$ small enough. For a introduction to the FvK-model and the definition of the energy $I_{h,\Delta}^{vK},$ see section 2.1. However a small correction needs to be made. The energy we give below in case of the excess cone differs from the one of the regular cone by a change of sign in front of the term related to the metric.\\
Let us note that Heiner Olbermann was able to improve this statement. It is optimal in the sense that the constants in front of the main term $h^2|\log h|$ agree in both bounds. In our case we will lose optimality. More details can be found in the remark below lemma \ref{lb22}. 
\\

\begin{figure}
    \centering
  \includegraphics[width=6cm,height=4cm]{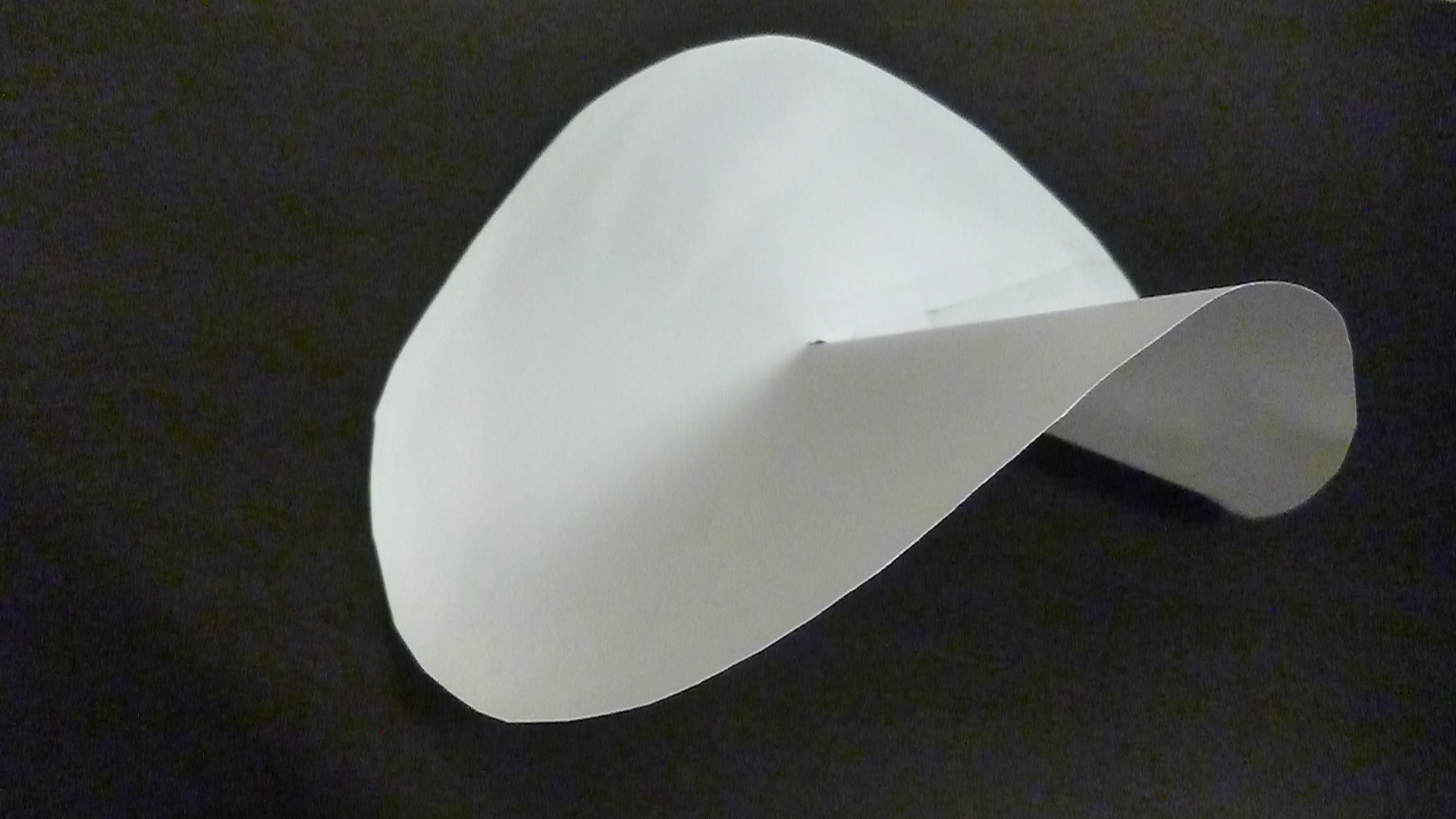}
    \caption{e-cone: additional angle approx. $45^\circ.$}
    \label{esmall}
\end{figure}

A regular cone can easily be constructed out of a piece of paper. We start by cutting out a circle. Then we cut out a circular sector of some angular. Gluing the edges of the resulting gap of the circle together leads to a regular cone. We assume that the glue doesn't change the surface or the material. Here we are interested in a very similar construction. Instead of cutting a piece out of the circle we add an additional sector. Back to the construction this means we take a second circle and cut straight from the boundary to the center and glue the circular sector we have left from the first construction in between the slit of the circle. This object has `too much angle' and starts to fold. Due to this, it is sometimes called an excess cone or `e-cone' for short, see figure \ref{esmall} for a small impression. On a technical level, the difference between the metric of the regular cone $g_c=Id-\Delta^2\xpp$ and the one of the e-cone $g_e=Id+\Delta^2\xpp$ is a sign change, see chapter 2 for notation.\\

Our aim is to prove a similar estimate as in (\ref{eq:OlB1}) in the case where the reference metric is that of an e-cone. For the main result see theorem \ref{UP1MIN}. For simplicity we use the Föppl-von Kármán(FvK) model. We start chapter 2 with an heuristic deduction of the FvK-model from the more general fully-non linear model. Afterwards we prove an upper bound with methods already used in \autocite{BKN13, MO14}. Since the regular cone and the cone are very closely related we can use the lower bound and the proof given in \autocite{Ol15} with slight modifications.\\

\begin{figure}
\centering
  \includegraphics[width=6cm,height=4cm]{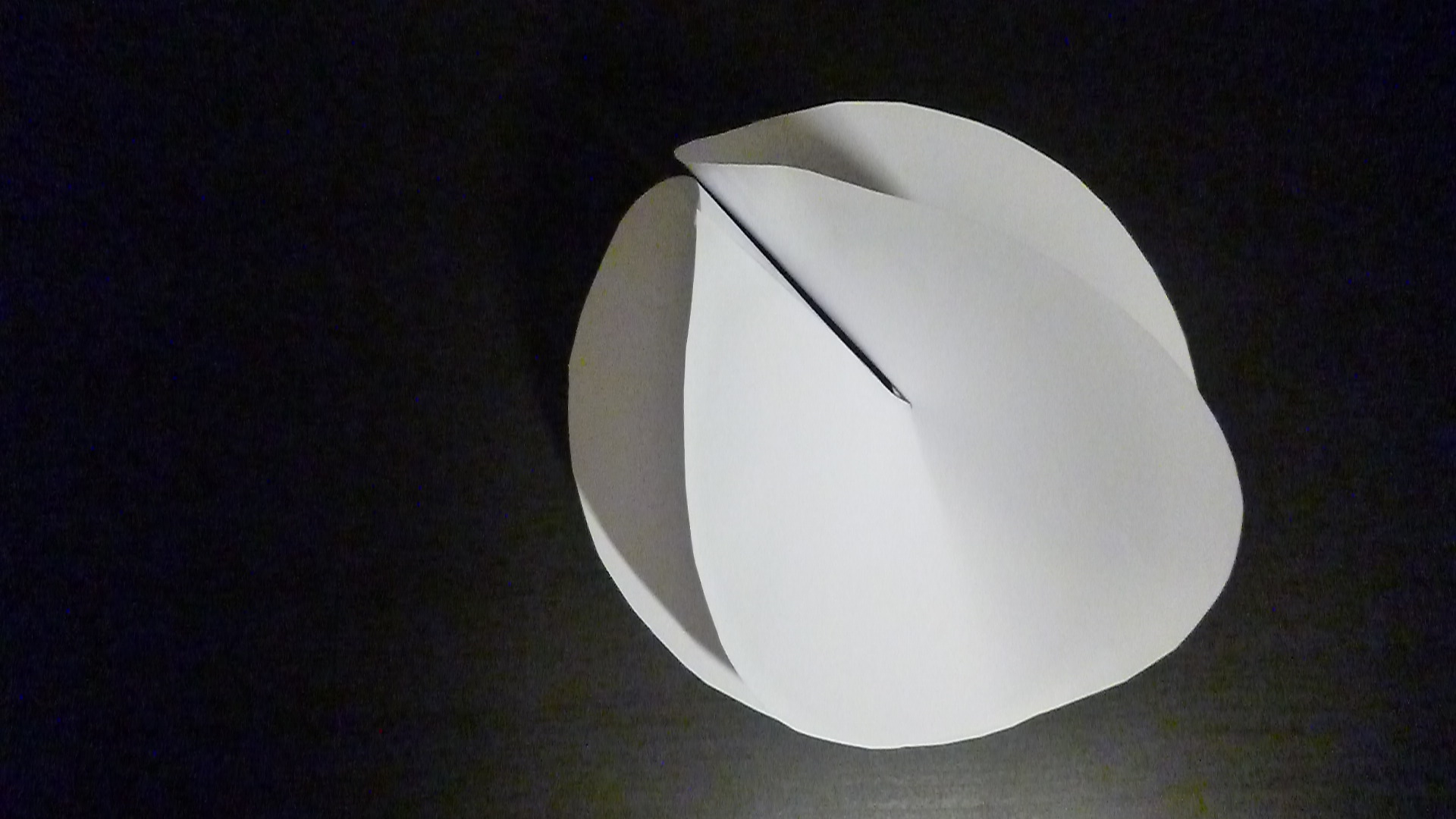}
    \caption{e-cone: additional angle approx. $360^\circ.$}
    \label{elarge}
\end{figure}

A theoretical and numerical analysis of the e-cone can be found in \autocite{MAG08}. They mainly considered one characteristic behavior of e-cones. With growing angle the folds become higher and closer together (see figure \ref{elarge}). At some point these folds will meet. The paper can't intersect itself. That's why the paper `reacts' by increasing the oscillation. This reduces the height of the folds. We will not see this phenomenon in the FvK-model, see the remark at the end of chapter 2.2.\\
\begin{figure}
    \centering
  \includegraphics[width=6cm,height=4cm]{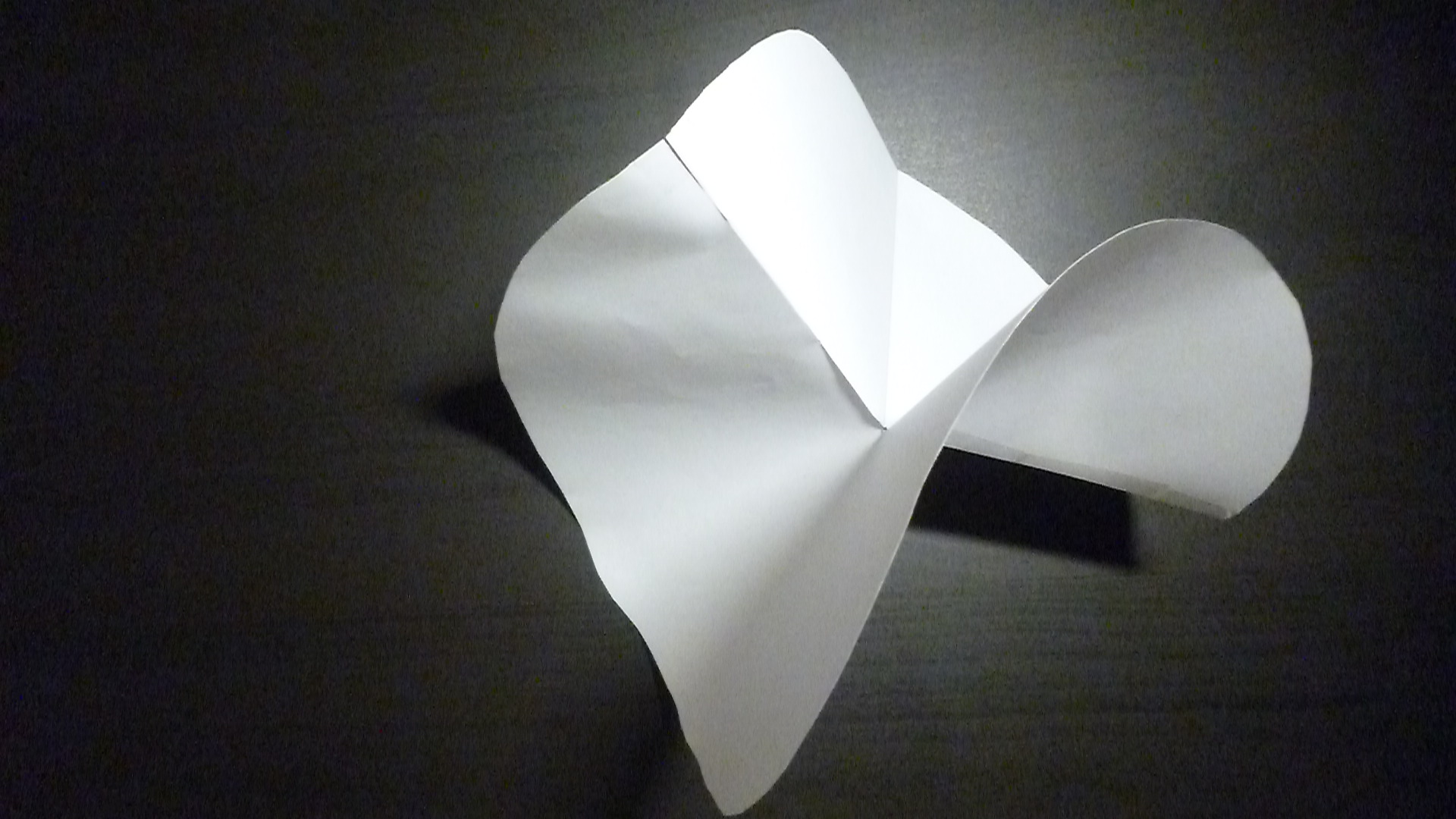}
    \caption{e-cone: additional angle approx. $90^\circ.$}
    \label{emiddle}
\end{figure}

 
In the second part we do not let the geometry evolve itself instead we interact with our object. For this imagine a spherical cap out of a thin elastic material. With the flat side facing down we position the cap on a table. Now we attach a sharp pencil pointing downwards on the highest point of the cap and start pushing perpendicular to the surface.\\

Now the questions are: how does the spherical surface behave under compression of the pencil and how does the energy scale, in dependence of a certain depth of the indentation $\delta$ and the thickness $h$? In the geometrically fully non-linear plate model or even in the 3D-model without any further restrictions these are very difficult questions which remain open. That is why we will restrict ourselves to the Föppl-von Kármán model and to radial symmetric functions.\\

In \autocite{COT15}, Conti, Olbermann and Tobasco considered the same restrictions but instead of looking at a spherical cap they considered a regular cone. They looked at the case where the cone was pushed at the top in the same manner as before. They showed that in case of radial symmetry the paper reacts in mainly three different ways depending on the depth of indentation $\delta$. If there is no push the energy scales like $h^2|\log h|$ and the cone remains the same. Considering a small push the cone starts bending on the base (outer boundary) and the energy behaves like $\delta^2h^{1/2}.$ For large indentations the cone inverts around the pencil and the corresponding energy scales with $\delta^{1/2}h^{3/2}.$\\
The upper bound is proven by the concrete construction of the above cases and estimation of the corresponding energy. To establish a lower bound it is crucial to understand the minimization of the energy as a two-well problem. This leads to different types of functions making it possible to show different kinds of lower bounds. We follow the ideas of the latter paper very closely.\\

In section $3$ we want to do the same approach in case of a spherical cap. As explained before in case of the regular cone we have three different regimes. However now in case of the cap we only have two. The characteristic energy of the unindented sphere scales like $h^2$ instead of the typical $h^2|\log h|$ term for cones. This time under compression however small spherical inversion occurs. The corresponding energy scales like $\delta^{3/2}h^{3/2}.$ For the concrete statement, see theorem \ref{THMSCmain}. \\
The regime for small indentations of the regular cone has vanished. This is surprising since this regime arises from the linear theory describing the small indentations as perturbations of the undeformed object. We would expect a similar reasoning for the spherical cap. However the cap seems very unstable, i.e. the non-linear term $\delta^{3/2}h^{3/2}$ is very small. Therefore we don't see the linear regime. On a local level the singularity leads to the rigidity of the singular cone.\\

There is another fascinating question raised and answered in \autocite{COT15}: if we drop the assumption of cylindrical symmetry, would we still see the three regimes (in case of the regular cone)? The answer is yes, for small indentations $0\le\delta\le h^{1/3}$. But under larger compression $\delta\ge h^{1/3}$ the symmetry gets broken. More precisely, it was shown that there are non-radial symmetrical constructions, which energy scales at most like $h^{5/3},$ which is smaller compared to the bound in the radial symmetric case for $\delta\ge h^{1/3}.$ We would expect the same behavior in case of the spherical cap (although we don't give a proof here).\\

In the same manner Ian Tobasco recently examined, in \autocite{T16}, the situation of a three dimensional thin elastic cylinder being compressed onto a rigid cylinder with same or larger radii, again under the assumption of cylindrical symmetry and showed energy scaling laws for each situation.\\

The deformation of spherical caps (or shells in $3D$) or other thin plates is not a newly arisen problem. It is of great interest to constructors, engineers in particular, aerospace engineers as well as physicists. For instance a theoretical discussion, closely related to ours, can be found in the book by Landau and Lifschitz \autocite{LL91} or in much greater detail in A.V. Pogorelov's book \autocite{Po88}. Pogorelov considers radial symmetrical deformations as we do, however his restrictions go beyond that. He restricts his argumentation to 'revolutions'(we call this spherical inversion) in order to be able to minimize the energy over the set of revolutions. Even if we know that such inversions are the minimizer, we allow any radial symmetric function. So, our approach is more general.\\
Pogorelov's analysis leads to the final energy $U$(with dimension $[U]=force\times length$)
\begin{equation}
U=2\pi c E\delta^{3/2}h^{5/2}/R.
\label{eq:}
\end{equation}
where $c$ is some (dimensionless) constant, $E$ Young's modulus with $[E]=pressure=force/(area)$, and $R$ the radius of the spherical shell, $[R]=length$. Of course $[\delta]=[h]=length.$\footnote{Note that in \autocite{Po88} $\delta$ and $h$ are used exactly the other way round, we changed this already.} Instead of looking at the energy we consider the energy per unit thickness $I,$ $[I]=force/h.$ Roughly, we claim that in our case $I=Ch^{3/2}\delta^{3/2},$ where $C$ is a constant with $[C]=force/(length)^3.$ If we now set $R=1\; length$ and $E=1\;force/(area)$ and we divide $U$ by $h\times length$ then the two results agree.\\

An experimental and numerical investigation of point indentations of spherical shells can be found in \autocite{NR14}. The results of Nasto and Reis confirm ours.\\
Other fascinating experiments have been made over the last few years concerning the deformation of thin shells (mainly metallic ones but also Ping Pong balls). To give a range of experiments, see \autocite{GPG99, GSV08, RGY06} or in \autocite{GSV07, DGY08} under dynamical compression. Since there are numerous different situations these experiments cannot be compared one by one to our setting. Yet most of them have the inversion behavior for small indentations in common and they all very clear show structures with broken symmetry under strong compression.

%% file: MT11.tex
\section{Energy scaling law for the singular e-cone}
In this first section we prove upper and lower bounds for the energy functional of the Föppl-von Karman picture. First we introduce the geometrically fully non-linear plate model together with some notation.\\

Let $B_1:=\{x\in\bb{R}^2:|x|<1\}$ be the unit ball with center $0$ in $\R^2$ and denote its boundary by $S^1.$   
For $x\in \R^2$ define  $\xh=\frac{x}{\xa},$ $\xp=\frac{1}{\xa}(-x_2,x_1).$ We will often use polar coordinates for this we denote the standard unit vectors of $\R^n$ by $e_1,\ldots,e_n.$ Further we denote polar coordinates in $\R^2$ by $e_r=(\cos t,\sin t)$ and $e_\vp=(\cos t,\sin t)$ for all $t\in[0,2\pi]$ and the cylindrical coordinates by $e_r,e_\vp,e_3$ respectively.

Then let further
\[g_y:=(D y)^TDy\]
be the metric induced by a map $y\in W^{2,2}(B_1,\bbs{3})$.\\

Fix some $h>0.$  For all $y\in W^{2,2}(B_1,\bbs{3})$ we define the energy functional 
\begin{eqnarray}\label{eq:GE}
I_h(y)=\int\limits_{B_1}{\left(|g_y-g_0|^2+h^2|D^2y|^2\right)\;dx.}
\end{eqnarray}
where we call $g_0$ the reference metric and the first term of the integral `membrane energy'. For it measures the deflection between the induced metric and the reference metric. Roughly speaking, the deflection of the two membranes. The second term including the second derivative and therefore measures the bending of the surface induced by $y$ and is therefore called `bending term'.\\

The above functional is indifferent under rotations. Indeed let $O(3)=\{R\in \R^{3\times3}:\; R^TR=Id_{3\times3}\}$ then it holds that 
\begin{equation}
I_h(Ry)=I_{h}(y)\; \mbox{for all}\; y\in W^{2,2}(B_1,\bbs{3})\;\mbox{and}\; R\in O(3).
\label{eq:FI}
\end{equation}
This can be seen by the small calculation  
\begin{eqnarray}
g_{(Ry)}&=&(DRy)^TD(Ry)= (D y)^TR^TRDy=(Dy)^TId_{3\times3} Dy=g_y\;\;\\
|D^2Ry|^2&=&|RD^2y|^2=R^TRD^2y:D^2y=Id_{3\times3}D^2y:D^2y=|D^2y|^2.\;\;
\label{eq:}
\end{eqnarray}
The reference metric $g_0$ does not depend on $y,$ hence the energy does not change under rotations. 
Geometrically this fact can be well understood. If we think of some object (like for example the e-cone below), the energy should not depend on how it is rotated in space. This property is crucial to the understanding of the result below.

In this section we are particularly interested in the case when the reference metric is that of a singular e-cone given by $g_e:=\Id+\Delta^2\xp\ot\xp.$ A very important class of functions will be the $1-$homogenous functions, i.e. functions of the form $\wt{y}(x)=|x|\gamma(\xh)$ for $\gamma\in C^2(S_1,\R^3).$ They are important because if $\gamma\in C^2(S_1,\R^3)$ satisfies the constraints
\begin{equation}
|\ga|=1\;\mbox{and}\;|\ga'|=\sqrt{1+\Delta^2}
\label{eq:BCFNP}
\end{equation}
then the induced metric of the corresponding $\wt{y}$ is that of a singular e-cone. In other words, $\wt{y}$ minimizes the membrane energy.
\subsection{Translation into the Föppl-von Kármán picture}
We don't want to use this general model, but instead we want to answer the question in the Föppl-von Kármán model. To see that the FvK model is of interest we at least want to argue heuristically why this model is important and how it could be related to the fully non-linear model.\\

For this sake we start from (\ref{eq:GE}) and we assume for small deflections, that $y\in C^{\infty}(B_1,\bbs{3})$ can be approximated up to 2nd-order by
\begin{eqnarray}\label{eq:}
y(x)=x+\ve^2U(x)+\ve W(x)e_3+O(\ve^3)
\end{eqnarray}
for small $\ve$ and $U\in C^\infty(B_1,\bbs{2})$ and $W\in C^\infty(B_1)$ where $U$ can be understood as the `inplane-' and $W$ as the `out-of-plane' deformation.

A small calculation shows for the first and second derivative of $y:$
\begin{eqnarray}\label{eq:}
\grad{y}(x)=\Id+\ve^2\grad{U}(x)+\ve e_3\ot\grad{W}(x)+O(\ve^3)\\
\grad^2{y}(x)=\ve^2\grad^2{U}(x)+\ve e_3\ot\grad^2{W}(x)+O(\ve^3)
\end{eqnarray}
Hence the bending term up to 2nd-order reads
\begin{equation}
|\grad^2y|^2=\ve^2 |\grad^2{W}|^2+O(\ve^3)
\label{eq:}
\end{equation}
and the induced metric becomes
\begin{equation}
g_y=\Id+\ve^2(\grad U+\grad^T U+\grad{W}\ot\grad{W})+O(\ve^3).
\label{eq:}
\end{equation}
Plugging this in the energy we get 
{\small
\begin{equation}\label{eq:FvKG}
I_h(y)=\ve^4\int\limits_{B_1}{\left(|2sym DU+\grad{W}\ot\grad{W}+\frac{1}{\ve^2}(\Id-g_0)|^2+h'^2|\grad^2W|^2\right)\;dx}+O(\ve^5)
\end{equation}}
where we redefined $h$ through $h':=\frac{h}{\ve}$. We want this term to be small but not so small that it is of order $\ve^5$ and therefore irrelevant, so assume $h'>\ve$. By $sym DU=(DU+DU^T)/2$ we denote the symmetrical part of the derivative of $U$.\\

Now we take as a reference metric, corresponding to a singular e-cone, given by $g_e:=\Id+\Delta^2\xp\ot\xp$ with some constant $0<\Delta<\infty,$ 
So the energy functional becomes
{\small
\begin{equation}\label{eq:}
I_h(y)=\ve^4\int\limits_{B_1}{\left(|2sym DU+\grad{W}\ot\grad{W}-\Delta'^2\xp\ot\xp|^2+h^2|\grad^2W|^2\right)\;dx}+O(\ve^5)
\end{equation}}
again, redefining $\Delta'=\frac{\Delta}{\ve}$ but $\Delta'$ still satisfies $0<\Delta'<\infty.$  The first non-trivial term in this energy is now the Föppl-von Kármán energy and therefore given by
{\small
 \begin{eqnarray}\label{eq:}
\IvK=\int\limits_{B_1}{\left(|2sym DU+\grad{W}\ot\grad{W}-\Delta'^2\xp\ot\xp|^2+h'^2|\grad^2W|^2\right)\;dx.}
\end{eqnarray}}
How do the boundary conditions for the $1-$homogenous functions translate?
For this assume $\wt{W}$ to be a $1-$homogenous function of the form
\begin{eqnarray}
\wt{W}(x)&=&|x|\alpha(\xh)
\end{eqnarray}
with $\alpha\in C^2(S_1,\R).$ Again we would like to have a pair of functions $(\wt{U},\wt{W})$ which represent the e-cone in the FvK-model, i.e. minimize the membrane energy. Since $\wt{W}$ is $1-$homogenous, $\wt{U}$ needs to be of the same form yielding
\begin{eqnarray}
\wt{U}(x)&=&|x|(u(\xh)e_r+v(\xh)e_\vp)
\label{eq:}
\end{eqnarray}
where $u,v\in C^2(B_1,\R).$\\
From now on we will switch between functions defined on the circle $f:S_1\rightarrow\R^n$ and functions defined on the interval $f:[0,2\pi]\rightarrow \R^n$ satisfying the boundary condition $f(0)=f(2\pi),$ without further notice.\\

Now the two conditions in (\ref{eq:BCFNP}) need to be satisfied up to 1st-order, or since we square them up to 2nd-order and we get
\begin{eqnarray}
|(1+\ve^2 u)e_r+\ve^2v e_\vp+\ve \alpha e_3|^2&=&1+O(\ve^4)\\
|\ve^2(u'-v)e_r+(1+\ve^2(u+v'))e_\vp+\ve \alpha'e_3|^2&=&1+\ve^2\Delta'^2+O(\ve^4)\;\;\;\;
\label{eq:}
\end{eqnarray}
Comparing components of the two sides and ignoring higher order terms leads to the conditions
\begin{eqnarray}
u=-\frac{\alpha^2}{2}\;\mbox{and}\;u+v'=\frac{1}{2}(\Delta'^2-\alpha'^2).
\label{eq:}
\end{eqnarray}
Combining them leads to the constraint
\begin{equation}
v'=\frac{1}{2}(\Delta'^2+\alpha^2-\alpha'^2).
\label{eq:}
\end{equation}
Since $v$ is $2\pi-$periodic, the integrated version becomes
\begin{equation}
\int_0^{2\pi}{(\al'^2-\al^2)\;dt}=2\pi\Delta'^2.
\label{eq:}
\end{equation} 
This motivates the definition of the following set
\begin{equation}
\B=\left\{\al\in C^2(S_1,\R):\int_0^{2\pi}{(\al'^2-\al^2)\;dt}=2\pi\Delta^2\right\}.
\label{eq:}
\end{equation}
Note that we dropped the prime on $\Delta,$ since from now on we will only deal with the FvK model so there should be no confusion. Further we will drop the prime on $h'$ to. 
After all this preliminaries we are finally able to state the first main result.
\begin{thm}\label{Thm1}
There exists a constant $C=C(\Delta)>0$ such that
{\small
\begin{equation}
2\pi\Delta^2h^2\left(\log{\frac{1}{h}}-2\log{\log {\frac{1}{h}}}-C\right)\le\inf\limits_{(U,W)} \IvK\le 6\pi \Delta^2 h^2\left(\log{\frac{1}{h}}+C\right)
\label{eq:}
\end{equation}}
for small enough $h,$ where $(U,W)\in W^{1,2}(B_1,\R^2)\times W^{2,2}(B_1).$
\end{thm}
The bounds are not optimal in the sense that the constants in front of the leading order term $h^2|\log h|$ differ by a factor of $3.$ Moreover new results in \autocite{Ol} give hope that the deviation from the leading order term is actually of order $h^2.$ For the reason why we lack of optimality and how optimality will imply the latter one, see remark below lemma \ref{lb22}.  
\subsection{Upper bound for the e-cone}
We start by proving an upper bound.
\begin{lem}\label{UP1}
There is a constant $C_2=C_2(\Delta)>0$ s.t.
\begin{equation}
\inf\limits_{(U,W)}\IvK\le 6\pi \Delta^2 h^2\log{\frac{1}{h}}+C_2h^2
\label{eq:}
\end{equation}
where $(U,W)\in W^{1,2}(B_1,\R^2)\times W^{2,2}(B_1).$
\end{lem}
{\bf{Proof.}}\\
Follows directly from the following two Lemmas \ref{UP1GEN} and \ref{UP1MIN}.\hfill$\square$\\

First we prove that there exists an upper bound in the first place. 

\begin{lem}\label{UP1GEN}
There are two constants $C_1=C_1(\Delta),C_2=C_2(\Delta)>0$ s.t.
\begin{equation}
\inf\limits_{(U,W)} \IvK\le C_1 h^2\log{\frac{1}{h}}+C_2h^2
\label{eq:}
\end{equation}
where $(U,W)\in W^{1,2}(B_1,\R^2)\times W^{2,2}(B_1)$ and 
\begin{equation}\label{}
C_1=\inf \limits_{\alpha\in \B}C_1(\alpha)=\frac{1}{\ln {2}}\inf \limits_{\wt{W}(x)=|x|\alpha(\xh),\atop \alpha\in \B}\int\limits_{B_1\setminus B_{1/2}}{|\grad^2 \wt{W}|^2\; dx}.
\end{equation}
\end{lem}
\vspace{0.5cm}

This form of $C_1(\alpha)$ has been proven by Brandman, Kohn and Nguyen \autocite{BKN13} and was also used in the setting of d-cones by Heiner Olbermann and Stefan Müller \autocite{MO14}. 
We still include the proof at this point.\\

{\bf{Proof.}}\\
Let $\eta\in C^2([0,\infty))$  be a truncation satisfying $\eta(t)=t,$ if $t\ge 1$ and $\eta(t)=0$ for all $t\in [0,\frac{1}{2}]$.
We want to choose the 1-homogenous functions representing the e-cone. Since they lack of regularity at the origin we use the truncation. Now choose  
\begin{eqnarray}
W_{h_*}(x)&=&h_*\eta\left(\frac{|x|}{h_*}\right)\alpha(\xh).\\
U_{h_*}(x)&=&h_*\eta\left(\frac{|x|}{h_*}\right)(u(\xh)e_r+v(\xh)e_\vp)
\end{eqnarray}
where $u=-\alpha^2/2,v'+u=(\Delta^2-\alpha'^2)/2.$\\
Then we split up the integral in a ball with radius $h$ and it's complement. In $B_{h}$ we have the estimates
$|\grad W_{h}|\le C$ and $|\grad^2 W_{h}|\le C/h.$ However in $B_1\setminus B_{h},$  $W_{h}$ and $U_{h}$ agree with the 1-homogenous functions defined above, i.e. $W_{h}=\wt{W}$ and $U_{h}=\wt{U}.$ Since $u$ and $v$ satisfy the constraints the membrane energy vanishes in $B_1\setminus B_{h}.$ Together this implies
\begin{eqnarray}\label{eq:}
I_{h,\Delta}^{\mbox{{vK}}}(U_{h},W_{h})\!\!\!&=&\!\!\!\int\limits_{B_{h}}{|2sym DU_{h}+\grad{W}_{h}\ot\grad{W}_{h}-\Delta^2\xp\ot\xp|^2\;dx}\;\;\;\;\;\;\\
&&+h^2\int\limits_{B_{h}}{|\grad^2W_{h}|^2\;dx}+h^2\int\limits_{B_1\setminus B_{h}}{|\grad^2W|^2\;dx}\\
&\le& Ch^2+Ch^2+h^2\int\limits_{B_1\setminus B_{h}}{|\grad^2W|^2\;dx}
\end{eqnarray}
For the last term we use a decomposition in dyadic annuli and the scaling property of $\wt{W}.$ For this sake, define $\Omega_n(x)=2^nW(2^{-n}x).$ Then $\Omega_n=W$ if $2^{-n}\ge h$ and by a change of variables:\\
\begin{equation}
\int\limits_{B_{2^{-n+1}}\setminus B_{2^{-n}}}{|\grad^2W|^2\;dx}=\int\limits_{B_1\setminus B_{\frac{1}{2}}}{|\grad^2\Omega_n|^2\;dx}=\int\limits_{B_1\setminus B_{1/2}}{|\grad^2 \wt{W}|^2\; dx}.
\end{equation}
Hence
\begin{equation}
h^2\int\limits_{B_1\setminus B_{h}}{|\grad^2W|^2\;dx}\le h^2\sum\limits_{n\ge0; 2^n\ge h}\int\limits_{B_{2^{-n+1}}\setminus B_{2^{-n}}}{|\grad^2W|^2\;dx}=C_1h^2\log\frac{1}{h}.
\end{equation}
This concludes the proof.
\hfill$\square$\\

Now we determine the concrete numerical value of $C_1.$

\begin{lem}\label{UP1MIN}
The set of minimizer $\mathcal{M}$ to the minimization problem
\begin{equation}\label{}
\inf \limits_{\wt{W}(x)=|x|\alpha(\xh), \atop \alpha\in \B}\frac{1}{\ln(2)}\int\limits_{B_1\setminus B_{1/2}}{|\grad^2 \wt{W}|^2\; dx}
\end{equation}
is given by
{\small\begin{eqnarray}
\mathcal{M}=\left\{ \al:S_1\rightarrow \R\left|\begin{array}{l}\alpha(t)=c_0\cos(t)+c_1\sin(t)+c_2\cos(2t)+c_3\sin(2t)\\
c_0,c_1,c_2,c_3\in\R\;\mbox{and}\;c_2^2+c_3^2=\frac{2}{3}\Delta^2\end{array}\right.\right\}.
\label{eq:}
\end{eqnarray}}
The minimal value is
\begin{equation}
6\pi \Delta^2.
\end{equation}
\end{lem}

{\bf{Proof.}}
We may use polar coordinates and consider functions $\wt{W}:[0,1]\times [0,2\pi]\rightarrow \R$ (Same notation as before, analogously for $\alpha$). From now on we want $\B$ to be the corresponding set of functions in polar coordinates, i.e. $\B$ consists of all functions $ \alpha\in C^2([0,2\pi])$ with periodic boundary conditions $\alpha(0)=\alpha(2\pi)$ and $\alpha$ satisfies the constraint from above. The corresponding minimization problem is given by

\begin{equation}\label{}
\inf \limits_{\wt{W}(r,t)=r\alpha(t), \atop \alpha\in \B} \frac{1}{\ln(2)}\int\limits_{B_1\setminus B_{1/2}}{|\grad^2 \wt{W}|^2\; dx}.
\end{equation}
Now we start by calculating the integral
\begin{eqnarray}
\frac{1}{\ln(2)}\int\limits_{B_1\setminus B_{1/2}}{|\grad^2 \wt{W}|^2\; dx}&=&\frac{1}{\ln(2)}\int\limits_{\frac{1}{2}}^1\int\limits_{0}^{2\pi}{(\alpha+\alpha'')^2\; \frac{1}{r}\; dt\;dr}\label{eq:MINI1}\\
&=&\int\limits_{0}^{2\pi}{(\alpha+\alpha'')^2\; dt}.
\label{eq:MINI2}
\end{eqnarray}
So the former two dimensional problem has reduced by $1-$ homogeneity of the functions $\wt{W}$ to a $1-$ dimensional minimization problem.

Since every minimizer of the last integral is a critical point of the functional including a Lagrange multiplier with constraint we will search for those. 
Recall that the constraint in an integrated form is given by 
\begin{equation}
\int\limits_0^{2\pi}{(\alpha'^2-\alpha^2)\;dt}=2\pi \Delta^2.
\label{eq:Con1}
\end{equation}
The expanded functional including the constraint then takes the form
\begin{eqnarray}
\int\limits_{0}^{2\pi}{(\alpha+\alpha'')^2+ \lambda(2\pi\Delta^2+\alpha^2-\alpha'^2)\;dt}.
\label{eq:}
\end{eqnarray}
We know that a minimizer to this problem exists and by regualarity therory we know that this function will be in $C^\infty.$
The corresponding Lagrangian $\La:[0,2\pi]\times\R^3\rightarrow\R$ is given by 
\begin{equation}
\La(t,\al,\al',\al'')=(\al+\al'')^2+\lambda(2\pi\Delta^2+\al^2-\al'^2).
\label{eq:}
\end{equation}
Then the Euler-Lagrange-equation looks like
\begin{eqnarray}
\partial_\al\mathcal{L}-\partial_t\partial_{\al'}\mathcal{L}+\partial^2_t\partial_{\al''}\mathcal{L}=2((1+\lambda)\al+(2+\lambda)\al''+\al^{(4)})=0
\label{eq:}
\end{eqnarray}
and finally
\begin{equation}
\alpha^{(4)}+(2+\lambda)\alpha''+(1+\lambda)\alpha=0.
\label{eq:ELE33}
\end{equation}
This equation has the general solution
\begin{equation}
\alpha(t)=c_0\cos(t)+c_1\sin(t)+c_2\cos(\sqrt{1+\lambda}t)+c_3\sin(\sqrt{1+\lambda}t)
\label{eq:GS1}
\end{equation}
for $\lambda\not=0$ and for $\lambda=0$
\begin{equation}
\alpha(t)=c_0\cos(t)+c_1\sin(t)+c_2t\cos(t)+c_3t\sin(t).
\label{eq:}
\end{equation}
The latter solution doesn't satisfy the constraint and boundary conditions, so there is no solution to \ref{eq:ELE33} satisfy the constraint in the case $\lambda=0$ (similarly for the case that $\lambda\le-1$).
for some $\lambda\not=0, \lambda>-1.$ Further the boundary conditions imply that $\sqrt{1+\lambda}$ should be a natural number greater than $1,$ i.e.
\begin{equation}
\sqrt{1+\lambda}=n \;\mbox{for some}\; n\in \N\setminus\{0,1\}.
\label{eq:Cond1}
\end{equation}
By plugging the general solution (\ref{eq:GS1}) into (\ref{eq:Con1}) and using relation (\ref{eq:Cond1}), for any such $n$ we get an equation for $c_{2},c_{3}.$ Since they depend on $n$ we will denote them by $c_{2,n},c_{3,n}$ and the relation is given by
\begin{eqnarray}
(c_{2,n}^2+c_{3,n}^2)(n^2-1)=2\Delta^2.
\label{eq:}
\end{eqnarray}
Let $n\in \N\setminus\{0,1\}$ and denote by $\al_n$ some solution of the form (\ref{eq:GS1}) with $c_{2},c_{3}$ replaced by $c_{2,n},c_{3,n}.$ With this knowledge we are able to calculate integral (\ref{eq:MINI2}):
\begin{eqnarray}
\int\limits_{0}^{2\pi}{(\alpha_n+\alpha_n'')^2\; dt}=\pi(c_{2,n}^2+c_{3,n}^2)(n^2-1)^2=2\pi\Delta^2(n^2-1).
\label{eq:}
\end{eqnarray}
The functions $\al_n$ are the possible candidates, we only need to minimize over $n$ to get the minimal value:

\begin{eqnarray}
\inf \limits_{\wt{W}(r,t)=r\alpha(t), \atop \alpha\in \B} \frac{1}{\ln(2)}\int\limits_{B_1\setminus B_{1/2}}{|\grad^2 \wt{W}|^2\; dx}&=&\min\limits_{n\ge 2}{\int\limits_{0}^{2\pi}{(\alpha_n+\alpha_n'')^2\; dt}}\\
&=&\min\limits_{n\ge 2}{2\pi\Delta^2(n^2-1)}\\
&=&6\pi\Delta^2.
\label{eq:}
\end{eqnarray}
The minimum is obtained for $n=2$ this yields the desired frequency of the minimizer.
\hfill$\square$\\

\textbf{Remark:} 1. At first this result looks kind of surprising. We would expect a solution of the form $c_2\cos(2t)+c_3\sin(2t)$ with some condition on $c_2,c_3$ at least for small $\Delta.$ But obviously we can add a function of the form $c_0\cos(t)+c_1\sin(t)$ for arbitrarily large amplitudes $c_0$ and $c_1.$ To understand this recall equation (\ref{eq:MINI1})-(\ref{eq:MINI2}) \\
\begin{eqnarray}
\frac{1}{\ln(2)}\int\limits_{B_1\setminus B_{1/2}}{|\grad^2 \wt{W}|^2\; dx}=\int\limits_{0}^{2\pi}{(\alpha+\alpha'')^2\; dt}.
\end{eqnarray}
Assume $\wt{W}(x)=|x|\alpha(\xh)$ is a linear function then the left side is zero\footnote{More generally this is true for affine functions but the only 1-homogenous affine functions are the linear ones.}. In this case the right hand side has to be zero, too, which is only possible if $\alpha(t)+\alpha''(t)=0$ for all $t\in[0,2\pi].$ The solution to this equation is exactly $c_0\cos(t)+c_1\sin(t)$ for all $c_0,c_1\in\R.$ Assume additionally, $\alpha$ is a minimizer of the above minimization problem. Then we can add a linear function to $\wt{W},$ i.e. $x\mapsto\wt{W}(x)+\nu\cdot x$ with $\nu\in \R^{2}$ without changing the value of the integral or equivalently, we can add $\alpha$ by a function of the form $c_0\cos(t)+c_1\sin(t),$ again without changing the integral. Hence this new function needs to be a minimizer, too.\\
We can see this from another perspective. Above we have discussed that the bending energy is invariant under the addition of a linear function to a 1-homogenous one. Assume that $\wt{U}$ is chosen s.t. the membrane energy vanishes for $(\wt{U},\wt{W}).$ Since $x\mapsto\wt{W}(x)+\nu\cdot x$ is again a 1-homogeneous function for all $\nu\in\R^2$ we can find a new $\wt{U}'$ s.t. the membrane energy vanishes, see the discussion above. This means the FvK model is invariant under adding a linear to a 1-homogeneous function.\\
This can be viewed as a remnant of the invariance of the fully non-linear plate model under 3-dimensional rotations. Note that the above needs not to be true for general functions $(U,W).$\\
 
2. In the introduction we have discussed that due to non-intersection of the paper for large $\Delta$ or to say it differently, for large additional angle, the paper reacts by increasing the frequency of oscillation which reduces the height of the object. We don't see this effect in the Föppl-von Kármán model or more precisely, the frequencies of the minimizer above don't depend on $\Delta$ (which we would expect to be the case in the fully non-linear model). Indeed, in the FvK-picture the sheet reacts only by increasing the amplitude of the oscillation.  


%% file: MT11LB.tex
\subsection{Lower bound for the e-cone}
For the upper bound we could use one type of possible functions and reduce the problem to a $1-$dimensional one, due to the fact that the minimum of the energy is always smaller or equal then the energy of an arbitrary configuration. Proofs of lower bounds (in similar settings) are much more involved. Our energy bound needs to be smaller than the energy for any functions.\footnote{Another reason why the lower bound is much more involved comes from the Nash-Kuiper Theorem which (roughly speaking) states, that there is a huge amount of configurations with arbitrarily small membrane energy. This would violate any kind of lower bound therefore we need to show that such functions would have high bending energies. For a detailed discussion, see \autocite{Ol15}.} \\
Before defining the determinant of the hessian in it's very weak form, we outline the proof. The proof closely follows \autocite{Ol15}. The next definition and lemma corresponds to definition $1$ and lemma $3$ in \autocite{Ol15}. The proof of the lower bound below is a combination of proposition $2$ and proposition $4$ in the same paper. A few things need to be slightly modified due to the different geometry. The lower bound is obtained by first showing that if for some pair of functions $(U,W)$ the upper energy bound is satisfied then the Gauss-curvature $\displaystyle\kvK(r)=\int\limits_{B_r}{\det D^2W\;dx}$ is close to the `curvature' of the singular e-cone. This can be made quantitative, we get an estimate on the $L_1-$ norm of $\kvK-\pi\Delta^2$ on some annulus. Then we can estimate the energy from below and by the latter estimate together with a dyadic decomposition we can conclude the proof.\\
\begin{de} Let $W\in W^{2,2}(B_1).$
By $\det D^2W:C_c^\infty(B_1)\rightarrow \R$ we denote the distribution
\begin{equation}
\vp\mapsto-\frac{1}{2}\int\limits_{B_1}{d\vp\wedge(W_{,1}dW_{,2}-W_{,2}dW_{,1})}.
\label{eq:}
\end{equation}
\end{de}
In the following we would like to consider such determinants also for 1-homogenous functions which are not in $W^{2,2}.$ For these functions we may interpret the determinant as a radon measure (same notation).
\begin{lem}Let $\wt{W}:B_1\rightarrow \R,$ $\wt{W}(x)=|x|\alpha(\xh)$ and $\alpha\in\B$ as before. Then 
\begin{equation}
\det D^2\wt{W}=-\Delta^2\pi\delta_0.
\label{eq:}
\end{equation}
where $\delta_0:\vp\rightarrow \vp(0)$ is the Dirac-distribution.
\end{lem}
{\bf{Proof.}}\\
For the proof we calculate the differential form $\wt{W}_{,1}d\wt{W}_{,2}-\wt{W}_{,2}d\wt{W}_{,1}$ in polar coordinates. Assume $r>0$ then
we can write $W_{,1}=W_{,r}\cos \vp+r^{-1}W_{,\vp}\sin\vp$ and $W_{,2}=W_{,r}\sin \vp-r^{-1}W_{,\vp}\cos\vp.$ The differential form in polar coordinates is given by
\begin{eqnarray}
W_{,1}dW_{,2}-W_{,2}dW_{,1}&=&\left(\frac{W_{,r}W_{,r\vp}}{r}-\frac{W_{,\vp}W_{,rr}}{r}-\frac{W_{,r}W_{,\vp}}{r^2}\right)dr\\
&&\!\!+\left(W_{,r}^2-\left(\frac{W_{,\vp}}{r}\right)^2+\frac{(W_{,r}W_{,\vp})_{,\vp}}{r}\right)d\vp.
\label{eq:Diffform}
\end{eqnarray}
Since $\wt{W}$ is $C^2$ away from the origin we can calculate the above expressions. This leads to 
\begin{equation}
\wt{W}_{,1}d\wt{W}_{,2}-\wt{W}_{,2}d\wt{W}_{,1}=(\alpha^2-\alpha'^2+\alpha'^2+\alpha\alpha'') d\vp,
\label{eq:}
\end{equation}
where the $dr-$term has vanished. Since $\alpha'^2+\alpha\alpha''=(\alpha\alpha')'$ and $\alpha$ is periodic we can neglect this part. Moreover $\alpha$ satisfies the constraint yielding
\begin{eqnarray}
\wt{W}_{,1}d\wt{W}_{,2}-\wt{W}_{,2}d\wt{W}_{,1}&=&-2\pi\Delta^2d\vp\\
&=&-\frac{2\pi\Delta^2}{|x|^2}(x_1dx_2-x_2dx_1)\\
&=&2\pi\Delta^2\delta_0dx_1\wedge dx_2.
\label{eq:}
\end{eqnarray}
This proves the lemma.
\hfill$\square$\\

\begin{lem}\label{lb22}
Assume for $(U,W)\in C^{2}(B_1,\R^3)$ the upper bound for the energy holds true, i.e. $\IvK\le6\pi\Delta^2h^2(\log{\frac{1}{h}}+C).$\\ 
Then $(U,W)$ satisfies
\begin{eqnarray}
2\pi \Delta^2 h^2\left(|\log{h}|-2\log{|\log {h}|}-C\right)\le\IvK, \label{eq:lb2}
\end{eqnarray}
for small enough $h.$
\end{lem}

{\bf{Proof.}}\\
\textbf{Step 1:} $L^1-$bound of the Gauss--curvature.\\
We show that for all $R\in(h,1]$ and $h$ small enough the following bound
\begin{eqnarray}
\|\kvK+\pi\Delta^2\|_{L^1(h,R)}\le C h^{1/2}R^{1/2}|\log{h}|
\label{eq:lb21}
\end{eqnarray}
holds.\\

Define $F'(r)=\kvK(r)+\pi\Delta^2.$ We want to apply the standard interpolation inequality given by
\begin{equation}
\|F'\|_{L^1(h,R)}\le\|F\|_{L^1(h,R)}^{1/2}\|F''\|_{L^1(h,R)}^{1/2}+\|F\|_{L^1(h,R)},
\label{eq:SIF}
\end{equation}
which is a special case of the Gagliardo-Nirenberg-inequality, we will see another version later. The general inequality can be found in \autocite[1.4.8, Thm 1]{Maz11} and the references given there.\\
For this purpose it is necessary to write $F$ in a convenient way and estimate $F$ mainly through the membrane energy which itself can be controlled from above. For this sake, we first write the membrane energy and the Gauss-curvature in polar coordinates this will give us an impression how to choose $F$\footnote{An easy way to obtain $F$ would be to integrate $F'.$ Then we could estimate $F$ in terms of the bending term and finally by the upper bound for the global energy leaving us with the slightly worse bound $CR|\log{h}|$ which is not good enough for our purpose.}.
The membrane energy $E_{mem}(U,W)$ takes the form
\begin{eqnarray}
\int\limits_0^1{r\;dr}\int\limits_0^{2\pi}{d\vp\;\left(|2U_{r,r}+W_{,r}^2|^2+|2r^{-1}(U_{\vp,\vp}+U_r)+(r^{-1}W_{,\vp})^2-\Delta^2|^2\right.}\\
\!\!\!\left.+|2U_{\vp,r}+r^{-1}(U_{r,\vp}-U_\vp+W_{,r}W_{,\vp})|^2\right)
\label{eq:}
\end{eqnarray}
The Gauss-curvature can be represented by
\begin{equation}
\kvK(r)=\int\limits_{B_r}{\det D^2W\;dx}=\frac{1}{2}\int\limits_{\partial B_r}{W_{,1}dW_{,2}-W_{,2}dW_{,1}}.
\label{eq:}
\end{equation}
We know this differential form already, recall again equation (\ref{eq:Diffform}). Keep in mind we want to choose $F$ such that it will be easy to compare it with the above quantities. We define two functions for $h\le s,r\le1:$
\begin{eqnarray}
F_1(s):=\int\limits_h^s{dr}\int\limits_0^{2\pi}{d\vp\;(2U_{r,r}+W_{,r}^2)},\\
F_2(r):=\int\limits_0^{2\pi}{d\vp\;(2U_r+\frac{W_{,\vp}^2}{r}-\Delta^2r)}.
\label{eq:}
\end{eqnarray}
Finally define $F(s):=\frac{1}{2}(F_1(s)-F_2(s)).$  A short calculation shows that the derivative of $F$ agrees with the definition above, indeed $F'(r)=\kvK(r)+\pi\Delta^2$ and moreover it holds
\begin{equation}
F''(r)=\int\limits_0^{2\pi}{r\;d\vp \det D^2W}.
\label{eq:}
\end{equation} 
This concludes the construction of $F.$\\

We now start by showing an $L^{\infty}-$ bound on $F_1.$ Therefore we first use the Cauchy-Schwarz inequality and then Jensen's inequality to obtain  
\begin{eqnarray}
|F_1(s)|&\le& C\left(\int\limits_0^1{r \;dr}\int\limits_0^{2\pi}{d\vp\;|2U_{r,r}+W_{,r}^2|^2}\right)^{1/2}\left(\int\limits_h^1{\frac{dr}{r}}\right)^{1/2}\\
&\le& CE_h^{1/2}|\log h|^{1/2}\\
&\le& Ch|\log h|
\label{eq:}
\end{eqnarray}
for all $s\in[h,1].$ Moreover this implies an $L^1-$bound on $F_1:$
\begin{eqnarray}
\|F_1\|_{L^1(h,R)}\le CRh|\log h|.
\label{eq:}
\end{eqnarray}
A similar argument works for $F_2$ as well. Note that $\int_0^{2\pi}{U_{\vp,\vp}\;d\vp}=0$ so we can add this term to $F_2$ without changing it and estimate as before
{\small
\begin{eqnarray}
\int\limits_h^1{|F_2(r)|\;dr}=\int\limits_h^1{r dr\left|\int\limits_0^{2\pi}{d\vp\;(2r^{-1}(U_{\vp,\vp}+U_r)+\frac{W_{,\vp}^2}{r^2}-\Delta^2)}\right|}\hspace{1cm}\\
\le C\left(\int\limits_0^1{r \;dr}\int\limits_0^{2\pi}{d\vp\;|2r^{-1}(U_{\vp,\vp}+U_r)+\frac{W_{,\vp}^2}{r^2}-\Delta^2|^2}\right)^{1/2}\left(\int\limits_h^1{r\;dr}\right)^{1/2}\\
\le CRE_h^{1/2}\le ChR|\log h|^{1/2}.\hspace{6.5cm}
\label{eq:}
\end{eqnarray}}
Together we get
\begin{equation}
\|F\|_{L^1(h,R)}\le CRh|\log h|.
\label{eq:}
\end{equation}
Now we turn to the estimate of $\|F''\|_{L^1(h,R)}:$ This follows directly by the estimate
\begin{eqnarray}
\|F''\|_{L^1(h,R)}&\le&\int\limits_h^1{r\;dr\int\limits_0^{2\pi}{d\vp\;|\det D^2W|}}\\
&\le&\int\limits_{B_1}{|D^2W|^2\;dx}\\
&\le& C|\log h|.
\label{eq:}
\end{eqnarray}
By the standard interpolation (\ref{eq:SIF}) the claim follows. \\

\textbf{Step 2:} Proof of (\ref{eq:lb2}).\\
Since
\begin{equation}
\IvK\ge h^2\int\limits_{B_R\setminus B_r}{|\grad^2W|^2\;dx}
\label{eq:}
\end{equation}
holds for arbitrary $0<r<R\le1,$ it is enough to show a lower bound of the bending energy on one specific annulus.  
A first and important step will be to estimate it mainly through the Gaussian-curvature from below.\\

For this recall that for every $v\in C^2(\overline{B_1})$ the following inequality holds
\begin{equation}
\int\limits_{\partial B_r}{|\grad^2v|\;d\mathcal{H}^1}\ge\left(4\pi\left|\;\int\limits_{B_r}det D^2v\;dx\right|\right)^{1/2}.
\label{eq:Iso}
\end{equation}
The major ingredients to prove this inequality is a combination of the Sobolev-inequality for functions with bounded variation and degree theory, details can be found in \autocite{Ol15}. Using Jensen's inequality and (\ref{eq:Iso}) we get
\begin{eqnarray}
\int\limits_{\partial B_r}{|\grad^2W|^2\;d\mathcal{H}^1}&\ge&2\pi r\left(\;\int\limits_{\partial B_r}{|\grad^2W|\;d\mathcal{H}^1}\right)^2\\
&\ge&\frac{2}{r}\left|\int\limits_{B_r}det D^2W\;dx\right|=\frac{2}{r}|\kvK|(r).
\label{eq:}
\end{eqnarray}
In preparation of the dyadic decomposition we set $h_0=5h|\log h|^2$ and choose $J\in \N$ such that $2^Jh_0\le1-5h\le 2^{J+1}h_0.$ Further we set $R_j:=2^jh_0$ for all $j=0,\ldots,J+1.$ Note that by definition $[R_0,R_J]\subseteq[5h,1-5h].$ Now we integrate the inequality achieved before over the interval $[5h,1-5h]$ w.r.t. $r:$
\begin{eqnarray}
\int\limits_{B_{5h}\setminus B_{1-5h}}{|\grad^2W|^2\;dx}&\ge&2\int\limits_{5h}^{1-5h}{\frac{1}{r}|\kvK|(r)\;dr} \label{eq:LB1Re}\\
&\ge&2\int\limits_{R_0}^{R_J}{\frac{1}{r}|\kvK|(r)\;dr}\\
&\ge&2\pi\Delta^2\int\limits_{R_0}^{R_J}{\frac{1}{r}\;dr}-2\int\limits_{R_0}^{R_J}{\frac{1}{r}|\kvK(r)+\pi\Delta^2|\;dr} \label{eq:LB2Ea}\\
&\ge&2\pi\Delta^2(|\log h|-2\log{|\log h|}-C)\\
&&-2\int\limits_{R_0}^{R_J}{\frac{1}{r}|\kvK(r)+\pi\Delta^2|\;dr}.
\end{eqnarray}
In the latter inequality we calculated the integral and used the definition of $R_0,\;R_J$ and $h_0.$ 
The rest-term integral can now be decomposed into dyadic intervals and by the $L^1-$ bound (\ref{eq:lb21}) we get the estimate
\begin{eqnarray}
\int\limits_{R_0}^{R_J}{\frac{1}{r}|\kvK(r)+\pi\Delta^2|\;dr}&=&\sum\limits_{j=1}^{J}{\int\limits_{R_{j-1}}^{R_j}{\frac{1}{r}|\kvK(r)+\pi\Delta^2|\;dr}}\\
&\le&\sum\limits_{j=1}^{J}{R_{j-1}^{-1}\|\kvK+\pi\Delta^2\|_{L^1(h,R_j)}}\\
&\le& Ch^{1/2}|\log{h}| \sum\limits_{j=1}^{J}{R_{j}^{-1/2}}\\
&\le& Ch^{1/2}|\log{h}|h_0^{-1/2} \frac{{2}^{-1/2}-{2}^{-(J+1)/2}}{1-{2}^{-1/2}}\\
&\le& C.
\label{eq:}
\end{eqnarray}
This finishes the proof.
\hfill$\square$\\

\textbf{Remark:} Parallel to this work, Heiner Olbermann improved the lower bound in case of the regular cone, in \autocite{Ol}. Indeed he showed that 
\begin{equation}
2\pi\Delta^2h^2(\log{\frac{1}{h}}-C)\le \min I_{h,\Delta}\le 2\pi\Delta^2h^2(\log{\frac{1}{h}}+C)
\label{eq:}
\end{equation}
where the second term in the lower bound does not show up any more. The proof used in \autocite{Ol} is much more elegant, since it is no longer necessary to perform a dyadic decomposition, it is enough to consider one ring with conveniently chosen radii. A similar estimate as (\ref{eq:LB1Re})-(\ref{eq:LB2Ea}) can be established and the corresponding negative rest part can be estimated by the membrane energy. Then one is able to get an estimate for the membrane energy from above and finally for the complete energy from below. Optimality, i.e. the fact that $2\pi\Delta^2$ shows up in the upper bound and in (\ref{eq:LB1Re}), too, is crucial for this argument to work.\\

We lack of optimality due to the fact that we used the inequality (\ref{eq:Iso}) (or more precisely: the Sobolev-inequality) for general functions $v\in C^2(\overline{B_1}).$ For the regular cone this is good enough, but to increase the lower bound of the e-cone one has to find a convenient restriction to a class of functions $X\subset C^2(\overline{B_1})$ s.t. one has the chance to improve the constant. In our opinion, to understand which class $X$ to choose and how it is related to functions satisfying just the upper bound (because that's the only thing we assume) would need tremendous insight. Still it is an interesting question, especially in the above context, since improving this constant would not only enhance the lower bound by itself but would also allow us to apply this new method and we could drop the $h^2\log|\log h|$ term.\\

The information above, is enough to show the complete result but we outline the proof and show how the pieces fit together.\\

{\bf{Proof of Theorem \ref{Thm1}.}} The direct method of calculus guarantees that the minimum is attained. The upper bound was proven in lemma \ref{UP1}. Moreover by the latter we also know that there is at least one pair of functions $(U,W)\in W^{1,2}(B_1,\R^2)\times W^{2,2}(B_1)$ such that $I_{h,\Delta}(U,W)\le 6\pi \Delta^2 h^2\log{\frac{1}{h}}+Ch^2,$ namely the construction given in the proof of the upper bound. Further $C^2(B_1,\R^3)$ is dense in $W^{1,2}(B_1,\R^2)\times W^{2,2}(B_1)$ and it holds that if $(U_\ve,W_\ve)\rightarrow(U,W)$ in $W^{1,2}(B_1,\R^2)\times W^{2,2}(B_1)$  for $\ve\rightarrow 0$ then $I_{h,\Delta}(U_\ve,W_\ve)\rightarrow I_{h,\Delta}(U,W).$\footnote{To see this convergence it is crucial to notice that $\|\grad W\ot\grad W\|_{L^2(B_1)}=\|\grad W\|_{L^4(B_1)}^2\le C\|\grad W\|_{W^{1,2}(B_1)}^2$ due to compact Sobolev-embedding $W^{1,2}\hookrightarrow L^4$ in $\R^2.$} Hence, we may assume $(U,W)\in C^2(B_1,\R^3)$ satisfies the upper bound. The lower bound then follows from lemma (\ref{lb22}). \hfill$\square$

%% file: MT21.tex
\section{Indentation of an elastic spherical cap}
In the last section we, were interested in how the geometry (or the metric) alone leads to the folding of a thin sheet and how the energy behaves in dependence of the thickness $h.$ As announced earlier we now want additionally to compress our object in this case a spherical cap. This compression will be described by boundary conditions depending on the depth of the indentation $\delta,$ see the Definition of $\mathcal{A}_\delta.$ Since $\delta$ is an important parameter, our scaling law will depend on $\delta$ and $h.$\\

\textbf{Notation:} From now on most of the time we are no longer interested in numerical constants. Therefore we use the notation '$f\lesssim g$' meaning there exists a constant $C>0,$ which is independent of any variable whatsoever, such that $f\le Cg.$\\

Again we work in the Föppl-von Kármán picture this time the reference metric being that of a spherical cap, which in the limit of small deflections can be approximated up to 2nd order by $g_{\textbf{\betteris}}=\Id+4\ve^2|x|^2\xhh+O(\ve^3).$ The FvK-energy then takes the form\footnote{compare with equation (\ref{eq:FvKG}).}
{\small
\begin{eqnarray}\label{eq:}
\IvK=\int\limits_{B_1}{\left(|2sym D U+\grad{W}\ot\grad{W}-4|x|^2\xhh|^2+h^2|\grad^2W|^2\right)\;dx.}
\end{eqnarray}}
We now only consider radially symmetric configurations given by 
\begin{equation}\label{eq:}
U(x)=\frac{1}{2}u(\xa)\xh \;\;\mbox{and}\;\; W(x)=w(\xa)
\end{equation}
where the pair of functions $(u,w)$ is taken from the set 
\begin{equation}\label{eq:}
\mathcal{A}_\delta=\{(u,w):(U,W)\in W^{1,2}(B_1;\R^2)\times W^{2,2}(B_1), w(0)=0, w(1)=1-\delta\}.
\end{equation}
This finally leads to the energy functional $E_h:\mathcal A_{\delta}\rightarrow \R,$
\begin{equation}
E_h(u,w)=\int\limits_0^1{\frac{u^2}{r}+r\left(u'+w'^2-4r^2\right)^2+h^2\left(rw''^2+\frac{w'^2}{r}\right)\;dr}.
\label{eq:}
\end{equation}
which plays the central role in this section. We sometimes will use $E_{mem}$ and $E_{bend},$ which correspond to the membrane part and the bending part of this functional. We now state the main Theorem of this section. 
\begin{thm}\label{THMSCmain}
There exists a numerical constant $C>0$ such that
\begin{equation}
\frac{1}{C}(h^2+\delta^{3/2}h^{3/2})\le\min\limits_{\mathcal{A}_\delta} E_h\le C (h^2+\delta^{3/2}h^{3/2}).
\label{eq:}
\end{equation}
for all $h\in(0,1/2]$, $\delta\in[0,1].$
\end{thm}

\subsection{Upper bound}
Again we start by proving an upper bound:
\begin{lem}\label{SCUP1}
\begin{equation}
\min\limits_{\mathcal{A}_\delta} E_h\lesssim h^2+\delta^{3/2}h^{3/2}
\label{eq:}
\end{equation}
for all $h\in(0,1/2]$, $0<\delta\le1.$
\end{lem}
Like above we can choose one particular pair of functions $(u,w)$ and estimate the energy. How should $(u,w)$ be chosen?\\
To get a first impression lets look at the case where we don't push at all. 
This corresponds to the case $\delta=0$ and we can just use the function $w(r)=r^2$ and $u=0$ representing the spherical cap in the Föppl-von Kármán picture and we can calculate the energy explicitly
\begin{equation}
E_h(u,w)=h^2\int\limits_0^1{\left(rw''^2+\frac{w'^2}{r}\right)\;dr}=4h^2.
\label{eq:}
\end{equation}
The above calculation is easy. It is even simpler since we were able to chose $w(r)=r^2$ due to the fact that $W(\cdot)=|\cdot|^2\in W^{2,2}(B_1)$ a necessary condition for $(u,w)\in\A_\delta.$ This is not true for the singular cone. The defining function for the singular cone is $w(r)=r$ yielding $W(\cdot)=|\cdot|\notin W^{2,2}(B_1)$ therefore one is compelled to truncate around the origin. We don't need a truncation this time leaving us in a slightly better position not just above but also in the proof below.\\
We have chosen $u$ such that especially the second membrane term vanishes even the whole membrane term disappears.\\

{\bf{Proof.}}\\
The idea behind our choice of the function $w$ is to mix the functions $w(r)=r^2$ and $w(r)=-r^2.$ To do this appropriate we need to connect these functions by a smooth function with properties s.t. the estimate of the energy stays simple enough. For this we need a function $w_0\in C^{\infty}([-1,1])$ with the following properties: Let $0<2l<R$ and $R^2+l^2=\frac{\delta}{2}$ then
\begin{eqnarray}
w_0(\pm 1)=0,\;w_0'(+1)=2(R+l),\;w_0'(-1)=-2(R-l)\;\mbox{and}
\label{eq:Sm1}
\end{eqnarray}
\begin{eqnarray}
\int\limits_{R-l}^{R+l}{4r^2-|w_0'(\frac{r-R}{l})|^2\;dr}=0.
\label{eq:Sm2}
\end{eqnarray}
For instance consider the function 
\begin{equation}
v(x)=
\begin{cases} 
-2R(x+1)-l(x^2-1)&-1<x\le x_0\\
2R(x-2x_0-1)+l(x^2-2x_0^2+1)&x_0<x\le 1\\
\end{cases}
\label{eq:}
\end{equation}
with $x_0= R/l+\sqrt{1+R^2/l^2}$ s.t. the boundary conditions are satisfied. 
Define now
\begin{equation}
w_1'(r)=
\begin{cases} 
-2r & 0\le r\le R-l\\ 
w_0'(\frac{r-R}{l})  & R-l \le r\le R+l \\ 
2r  & R+l\le r\le 1
\end{cases}
\label{eq:}
\end{equation}
and we choose $u'_1$ like described before in such away that the `complicated' term $(u'+w'^2-4r^2)^2$ of the membran energy vanishes:
\begin{equation}
u_1'(r)=
\begin{cases} 
0 & 0\le r\le R-l\\ 
4r^2-w_0'^2(\frac{r-R}{l})  & R-l \le r\le R+l \\ 
0  & R+l\le r\le 1
\end{cases}
\label{eq:}
\end{equation}
and by integration we also have $w_1(r)=\int_0^r{w_1(t)\;dt}$ and
\begin{equation}
u_1(r)=\begin{cases} 
0 & 0\le r\le R-l\\ 
\int_{R-l}^r{(4t^2-|w_0'(\frac{t-R}{l})|^2)\;dt}  & R-l \le r\le R+l \\ 
0  & R+l\le r\le 1
\end{cases}
\label{eq:}
\end{equation}
Then it holds that $w_1(0)=0$ and
\begin{eqnarray}\label{}
w_1(1)&=&\int\limits_{0}^{R-l}-2r{\;dr}+\int\limits_{R-l}^{R+l}w_0'(\frac{r-R}{l}){\;dr}+\int\limits_{R+l}^{1}{2r\;dr}\\
&=&-(R-l)^2+l(w_0(1)-w_0(-1))+(1-(R+l)^2)\\
&=&1-2(R^2+l^2)\\
&=&1-\delta.
\end{eqnarray}
This yields $(u_1,w_1)\in \mathcal{A}_\delta,$ by our choice of $R$ and $l.$\\
Now we can establish estimates on each term of the energy separately. We start with the membrane terms, recall that we have chosen $u'$ in such a way that the second term vanishes. For the other one we first need a bound on $u_1$ which for $r\in[R-l,R+l]$ is given by
\begin{eqnarray}
|u_1(r)|&=&\left|\int_{R-l}^r{(4t^2-|w_0'(\frac{t-R}{l})|^2)\;dt}\right|\\
&\le&\frac{4}{3}((R+l)^3-(R-l)^3)+\int_{R-l}^{R+l}{|w_0'(\frac{t-R}{l})|^2\;dt}\\
&=&\frac{16}{3}(l^3+3R^2l)\\
&\lesssim&R^2l
\label{eq:}
\end{eqnarray}
where we first applied the triangle inequality and then used  (\ref{eq:Sm2}).
and therefore
\begin{eqnarray}
\int\limits_0^1{\frac{u_1^2}{r}\;dr}=\int\limits_{R-l}^{R+l}{\frac{u_1^2}{r}\;dr}
\lesssim R^4l^2\int\limits_{R-l}^{R+l}{\frac{1}{r}\;dr}
\lesssim R^3l^3
\label{eq:}
\end{eqnarray}
where we used 
\begin{equation}
\int\limits_{R-l}^{R+l}{\frac{1}{r}\;dr}\lesssim \frac{l}{R-l}\lesssim \frac{l}{R}
\label{eq:}
\end{equation}
by the relation $l<\frac{R}{2}.$\\

Now we consider the bending terms. Again using the $L^2-$property of the connection (\ref{eq:Sm2}) allows for the estimate 
\begin{eqnarray}
\int\limits_0^1{\frac{w_1'^2}{r}\;dr}&=&\int\limits_0^{R-l}{4r\;dr}+\int\limits_{R-l}^{R+l}{\frac{\left|w_0'(\frac{r-R}{l})\right|^2}{r}\;dr}+\int\limits_{R+l}^1{4r\;dr}\\
&\le&2(R-l)^2+\frac{1}{(R-l)}\int\limits_{R-l}^{R+l}{\left|w_0'\left(\frac{r-R}{l}\right)\right|^2\;dr}\\
&&+2(1-(R+l)^2)\\
&\le& 2+8Rl+\frac{16}{3}\frac{l^3}{R}\\
&\lesssim& 1.
\label{eq:}
\end{eqnarray}
On the other hand we have
\begin{eqnarray}
\int\limits_0^1{rw_1''^2\;dr}&=&\int\limits_0^{R-l}{4r\;dr}+\frac{1}{l^2}\int\limits_{R-l}^{R+l}{rw_0''^2\left(\frac{r-R}{l}\right)\;dr}+\int\limits_{R+l}^1{4r\;dr}\\
&\lesssim& 1+Rl+\frac{R^3}{l}.
\label{eq:}
\end{eqnarray}
Together this leads to the following bound on the energy
\begin{eqnarray}
E_h(u_1,w_1)&\lesssim& R^3l^3+ h^2(1+\frac{R^3}{l})
\label{eq:}
\end{eqnarray}
We don't want our bound in terms of $R$ and $l$ but in dependence of $h$ and $\delta$ only. For this sake, we have to minimize in $R$ and $l.$ 
If $\delta>h$ set $R=2^{-1/2}\sqrt{\delta-2^{-1}h}\le2^{-1/2}\sqrt{\delta},$ and $l=\frac{1}{2}h^{1/2}$ and we get
\begin{eqnarray}
E_h(u_2,w_2)&\lesssim&h^2+\delta^{3/2}h^{3/2},
\label{eq:}
\end{eqnarray}
if $\delta<h$ set $R=\frac{4}{5}2^{-1/2}\sqrt{\delta},\; l=\frac{3}{5}2^{-1/2}\sqrt{\delta}$ and we get
\begin{eqnarray}
E_h(u_2,w_2)&\lesssim&h^2+h^2\delta+\delta^{3}\lesssim h^2.
\label{eq:}
\end{eqnarray}
Hence for all $\delta\in[0,1]$ the upper bound is given by 
\begin{eqnarray}
E_h(u_2,w_2)&\lesssim&h^2+\delta^{3/2}h^{3/2}.
\label{eq:}
\end{eqnarray}
\hfill$\square$\\


%% file: MT22.tex
\subsection{Lower bounds}
To establish lower bounds we split the functions up into classes with different properties. An idea what kind of properties could be the right choice one can get by looking at the energy functional from the following point of view:

Assume $E_h(u,w)$ is small then $u$ is small in an $L^2-$ sense. And $u'$ is small in a weak sense. This means in this situation we can neglect all $u-$terms for a moment. Then the functional becomes 
\begin{equation}
\int\limits_0^1{(w'^2-4r^2)^2+h^2|w''|^2\;dr}.
\label{eq:}
\end{equation}   
So in a way, our problem is a version of the well known two-well-problem. This problem would be solved by the two solutions $w'(r)=\pm2r$. Of course we have ignored $u,$ so this will not exactly be the solution, but there is hope that 'good' neighborhoods of these solutions (we will refer to them as 'wells') will be of particular interest.  

We now introduce the same setting as in \autocite{COT15} but with slight changes of the definition of the wells and the function $\tau_f$ due to the differences in the energy functional. Note that this section including the lemmas follows this paper very closely.
 
The wells are now defined as follows for all $t\in [0,1]$ we set 
\begin{equation}
W_t=W_t^+\cup W_t^-,\;W_t^+=\left(\frac{3}{2}t,\frac{5}{2}t\right),\;\mbox{and}\;W_t^-=\left(-\frac{5}{2}t,-\frac{3}{2}t\right).
\label{eq:}
\end{equation}

For a function $f\in C([0,1])$ we denote by $\tau_f$ the last point before the function $f$ always remains in the wells, i.e. 
\begin{equation}
\tau_f=\max\{t\in(0,1]:f(t)\notin W_t\}.
\label{eq:}
\end{equation} 
If $f(t)\in W_t$ for all $t\in(0,1]$ then $\tau_f=0.$
We are interested in the case $f=w'$ for $w\in W^{1,2}([0,1])$ (More precisely, the continuous representative of $w'.$ That such a representative always exists, see  . We still stick to our notation) and set
\begin{equation}
\tau=\tau_{w'}.
\label{eq:}
\end{equation}
For the whole chapter we make the following assumptions 
\begin{equation}
\delta\in[0,1],\;h\in(0,1/2], \;(u,w)\in\A_\delta,\; E_h=E_h(u,w).
\label{eq:SetLB2}
\end{equation}

\textbf{1. Full spherical inversion.} Since we are not able to proof a lower bound for all functions at once we will show different estimates for different categories of functions. We first turn to the class of functions, which are going into the wells for the last time, on a very late time. We don't need to assume this for the next results but it is good to keep this in mind, since this is how we want to use it later.

Before proving the main bound in this case, we need to introduce a tool used throughout the proof, first. The next definition and lemma corresponds to \autocite[Lemma 8]{COT15}. There is a slight correction in the definition of the function $g_a$ due to the different geometry. This does not affect the proof but we will still include it here.
\begin{lem}\label{LB2:1}
Assume (\ref{eq:SetLB2}) holds. Let $w\in \WSob$ and $a\in[0,\frac{1}{2}].$ Define
\begin{equation}
I_a=[a,2a],\; g_a(r)=\int\limits_a^r{(4t^2-w'^2)\;dt}+c_a
\label{eq:}
\end{equation}
and $c_a$ chosen in such a way that the mean value of $g_a$ vanishes, i.e. $\fint\limits_{I_a}{g_a\;dt}=0.$
Then
\begin{equation}
\|g_a\|_{L^2(I_a)}\lesssim a^{1/2}E_h^{1/2}.
\label{eq:}
\end{equation}
\end{lem}
{\bf{Proof.}}\\
Let $v=u-g_a$ and define $v_a=\fint\limits_{I_a}{v\;dr}$ 
\begin{eqnarray}
E_h\ge E_{mem}&\ge&\int\limits_{I_a}{\frac{(v+g_a)^2}{r}+r|v'|^2\;dr}\\
&\ge&\frac{1}{a}\int\limits_{I_a}{(v+g_a)^2\;dr}+a\int\limits_{I_a}{|v'|^2\;dr}\\
&\ge&\frac{1}{a}\int\limits_{I_a}{(v+g_a)^2+|v-v_a|^2\;dr}
\label{eq:}
\end{eqnarray}
where we used Poincaré's inequality with optimal constant $\|v-v_a\|_{L^2(I_a)}\le a\|v'\|_{L^2(I_a)}.$ Expanding the latter term and using Cauchy's inequality for the mixed terms 
\begin{equation}
2(g_a-v_a)v\ge-\frac{1}{2}(g_a-v_a)^2-2v^2\ge-\frac{1}{2}g_a^2-\frac{1}{2}v_a^2+g_av_a-2v^2.
\label{eq:}
\end{equation} 
Note that the term $g_av_a$ has mean zero, hence
\begin{eqnarray}
E_h\ge\frac{1}{2a}\int\limits_{I_a}{g_a^2+v_a^2\;dr}\gs\frac{1}{a}\int\limits_{I_a}{g_a^2\;dr}
\label{eq:}
\end{eqnarray}
which completes the proof.\hfill$\square$\\
We now proof a lower bound of the energy in terms of $\ta.$
\begin{lem}\label{LB2:2}
Assume (\ref{eq:SetLB2}) holds and $\ta$ defined as above for some $w\in \WSob$. Then
\begin{equation}
\min\{\tau^6,\tau^{3}h^{3/2}\}\ls E_h.
\label{eq:}
\end{equation}
\end{lem}
{\bf{Proof.}}\\
We may assume $\tau>0.$ (If $\tau=0$ the statement $E_h\ge0$ is always true.) By definition of $\tau,$ $f(\tau)\notin W_\tau$ yielding the bounds $|w'(\ta)|\le\frac{3}{2}\tau$ or $|w'(\ta)|\ge\frac{5}{2}\tau.$
Let $0<\ve\le \frac{\tau}{2}$ which will be chosen explicitly below, and define $\Ite=[\ta-\ve,\ta]\subseteq[\frac{\ta}{2},\ta].$ Additionally we assume
\begin{equation}
\osc{\Ite}{w'}\le \frac{1}{4}\tau
\label{eq:OSC}
\end{equation}
this implies $|w'|\le\frac{7}{4}\ta$ or $|w'|\ge\frac{9}{4} \ta$ on $\Ite.$  Now by the latter inequalities we get $|g_{\frac{\ta}{2}}'|\gs\ta^2$ on $\Ite,$ where $g_{\frac{\ta}{2}}'$ does not change sign due to continuity. Integration leads to the scaling $|g_{\frac{\ta}{2}}|\gs\ta^2\ve$ on $\Ite.$ With this at hand we are able to estimate the energy from below:
\begin{equation}
E_h\gs\frac{1}{\ta}\int\limits_{\frac{\ta}{2}}^{\ta}{g_{\frac{\ta}{2}}^2\;dr}\gs\frac{1}{\ta}\int\limits_{\Ite}{g_{\frac{\ta}{2}}^2\;dr}\gs\ve^3\ta^3.
\label{eq:}
\end{equation}
On the other hand if (\ref{eq:OSC}) does not hold we can estimate 
\begin{equation}
\frac{E_h}{h^2}\gs\int\limits_0^1{rw''^2\;dr}\gs\int\limits_{\Ite}{rw''^2\;dr}\gs\tau\int\limits_{\Ite}{w''^2\;dr}\gs\frac{\ta}{\ve}\left(\osc{\Ite}{w'}\right)^2\gs \frac{\ta^3}{\ve},
\label{eq:}
\end{equation}
hence
\begin{equation}
E_h\gs\frac{\ta^3h^2}{\ve}.
\label{eq:}
\end{equation}
Together with the bound achieved first we get
\begin{equation}
E_h\gs\min\{\ve^3\ta^3,\frac{\ta^3h^2}{\ve}\}.
\label{eq:}
\end{equation}
Setting $\ve=\min\{h^{1/2},\tau/2\}$ yields
\begin{equation}
E_h\gs\min\{\tau^6,\tau^{3}h^{3/2}\}.
\label{eq:}
\end{equation}
\hfill$\square$\\
Assume the energy is very small. Then by the statement above $\tau$ has to be small. Low energy forces the configuration to remain in the wells nearly the whole time.\\

\textbf{2. Spherical lower bound.} In section 3.1 we have already seen that for $\delta=0$ the unindented sphere satisfies $E_h\ls h^2.$ Recall the lower bound we try to achieve $h^2+\delta^{3/2}h^{3/2}.$ If $\delta$ is close to zero then the $h^2$ term will be the leading term. Therefore we need to show a bound of the form $h^2\ls E_h.$\footnote{Note that as a result this shows an energy scaling law for the (undeformed) spherical cap under radial symmetry of the form $C'\le \frac{E_h}{h^2}\le C$ where $0<C'\le C$ are two constants.} This is the statement of the next lemma. 
\begin{lem}\label{LB2:SpB}
Assume (\ref{eq:SetLB2}). Then the following bound holds
\begin{eqnarray}
h^2&\ls&E_h.
\label{eq:}
\end{eqnarray}
\end{lem}
{\bf{Proof.}}\\
Let $\ta$ be as before, and set $\ta_0=\frac{1}{2}$. If $\ta<\ta_0$ then $w'(r)\in W_r$ for all $r\in[\ta_0,1].$ Then we can estimate the energy from below by
\begin{equation}
\frac{E_h}{h^2}\ge\int\limits_{\ta_0}^{1}{\frac{w'^2}{r}\;dr}\gs\ta_0^2\log\frac{1}{\ta_0}\gs 1.
\label{eq:}
\end{equation} 
For the other case $\ta\ge\ta_0$ using Lemma (\ref{LB2:2}) we obtain
\begin{equation}
E_h\gs\min\{\tau_0^6,\tau_0^{3}h^{3/2}\}\gs\min\{1,h^{3/2}\}= h^{3/2}\gs h^{2}.
\label{eq:}
\end{equation}
\hfill$\square$\\

\textbf{2. Partially inverted lower bound.} For the next class of functions consider again our construction from the upper bound. This function starts for small $t$ in the well $W_t^-.$ Then it jumps around $R$ from the well $W_t^-$ to $W_t^+$ and remains in this well up to the boundary. Since $R\ls\sqrt{\delta}$ is small, the configuration remains in $W_t^+$ a very long time.\\
We are now interested in functions which stay in the wells from a very early (in some sense) time on.\\ Before that we have to prove a few inequalities. This corresponds to lemma 10 in \autocite{COT15}. The proof remains the same but there are slight changes due to the different metric.
\begin{lem}\label{LB2:3}
Assume (\ref{eq:SetLB2}) holds and $\ta$ defined as above for some $w\in \WSob$. Further let $0<a<\frac{1}{2}$ and $w'(r)\in W_r$ for all $r\in I_a.$ Then the following bounds hold
\begin{eqnarray}
\int\limits_{I_a}{|4t^2-w'^2(t)|^2\;dt}&\ls&\left(\frac{1}{a}+\frac{a}{h}\right)E_h, \label{eq:LB2:3.1}\\
\osc{I_a}{g_a}&\ls&\left(1+\frac{a^{1/2}}{h^{1/4}}\right)E_h^{1/2}.
\label{eq:LB2:3.2}
\end{eqnarray}
\end{lem}
{\bf{Proof.}}\\
Recall another special case of the Gagliardo-Nirenberg-interpolation inequality:
\begin{eqnarray}
\|g_a'\|_{L^2(I_a)}\ls\|g_a\|_{L^2(I_a)}^{1/2}\|g_a''\|_{L^2(I_a)}^{1/2}+\frac{1}{a}\|g_a\|_{L^2(I_a)}.
\label{eq:}
\end{eqnarray} 
By the triangle inequality and since $w'(r)\in W_r$  for all $r\in I_a$ we get
\begin{eqnarray}
\|g_a''\|_{L^2(I_a)}^2&=&\int\limits_{I_a}{|8t-2w'(t)w''(t)|^2\;dt}\\
&\ls&\int\limits_{I_a}{t^2\;dt} +\int\limits_{I_a}{t^2|w''(t)|^2\;dt}\\
&\ls&a^3 +a\int\limits_{I_a}{t|w''(t)|^2\;dt}\ls a^3+\frac{aE_h}{h^2}\ls\frac{aE_h}{h^2}
\label{eq:}
\end{eqnarray}
where we used $a^2\ls1\ls \frac{E_h}{h^2}$ in the last step which is true by Lemma \ref{LB2:SpB}.
Finally we use Lemma \ref{LB2:1} to estimate $\|g_a\|_{L^2(I_a)}$ which implies the desired bound.

A similar argument works for the second inequality we just have to replace the above version of the Gagliardo-Nirenberg-inequality 
by
\begin{eqnarray}
\osc{I_a}{g_a}&\ls&\|g_a\|_{L^2(I_a)}^{3/4}\|g_a''\|_{L^2(I_a)}^{1/4}+\frac{1}{a^{1/2}}\|g_a\|_{L^2(I_a)}
\label{eq:}
\end{eqnarray} 
for details, see \autocite[Lemma 20]{COT15}. Using the same estimates as before get the claim. \hfill$\square$\\

The next lemma contains estimates on the $L_1-$norm of $w'$, see \autocite[Lemma 11]{COT15}. Like above the result remains nearly the same.
\begin{lem}\label{LB2:4}
Assume (\ref{eq:SetLB2}) holds and let $0<s\le t\le\frac{1}{2}.$ Then the following bounds hold
\begin{eqnarray}
\|w'\|_{L^1[0,s]}&\le& \frac{s}{h}E_h^{1/2},\\
\|w'\|_{L^1[s,t]}&\le& Ct^{1/2}(1+\log\frac{1}{s})^{1/4}E_h^{1/4}+\frac{8}{\sqrt{3}}t^2.
\label{eq:}
\end{eqnarray}
\end{lem}
{\bf{Proof.}}\\
The first estimate is just an application of Cauchy-Schwarz
\begin{eqnarray}
\int\limits_0^s{|w'|\;dr}\le\left(\int\limits_0^s{\frac{|w'|^2}{r}\;dr}\right)^{1/2}\left(\int\limits_0^s{r\;dr}\right)^{1/2}\le\frac{s}{h}E_h^{1/2}.
\label{eq:}
\end{eqnarray}
Now for the second inequality we need to be a little bit more careful. Before we have used the bending term and the relation $E_{bend}\le \frac{E_h}{h^2}$ now we use the membrane energy and the relation $E_{mem}\le E_h$ where there is no direct dependence on $h$ any more (of course $E_h$ still depends on $h$). Since the membrane term includes terms of $u$ and $u'$ we first need some estimates on them.
First note that by Tchebyshev's inequality for all $\ve>0$ it holds that
\begin{equation}
\La^1(\{r\in[0,1]:|u(r)|>\sqrt{\ve r}\})<\frac{1}{\ve}\int\limits_{0}^{1}{\frac{|u|^2}{r}\;dr}\le\frac{1}{\ve}E_h.
\label{eq:}
\end{equation}
To get the impression of a parameter we can choose arbitrarily, introduce $\gamma=\ve^{-1}E_h$ and the estimate above becomes
\begin{equation}
\La^1(\{r\in[0,1]:|u(r)|>(\frac{r}{\gamma})^{1/2}E_h^{1/2}\})<\gamma.
\label{eq:}
\end{equation}
Now choose $\gamma=\frac{s}{2}$ and $\gamma=t.$ Then there exists $r_0\in[\frac{s}{2},s]$ and $r_1\in[t,2t]$ s.t. the inequality from above is not true at these points and
\begin{equation}
|u(r_0)|+|u(r_1)|\le 3E_h.
\label{eq:}
\end{equation}                                      
Now we can start by applying the Cauchy-Schwarz inequality   
\begin{eqnarray}
\int\limits_s^t{|w'|\;dr}\le\sqrt{2t}\left(\;\int\limits_{r_0}^{r_1}{|w'|^2\;dr}\right)^{1/2}.
\label{eq:LB2CS}
\end{eqnarray}
The triangle inequality and Cauchy-Schwarz (again) lead to
\begin{eqnarray}
\!\!\!\!\!\!\!\!\int\limits_{r_0}^{r_1}{|w'|^2\;dr}\!\!\!\!&\le&\!\!\!\!\int\limits_{r_0}^{r_1}{u'+w'^2-4r^2\;dr}+\frac{4}{3}|r_1^3-r_0^3|+|u(r_1)-u(r_0)|\\
&\le&\!\!\!\!\left(\int\limits_{r_0}^{r_1}{r|u'+w'^2-4r^2|^2\;dr}\right)^{\frac{1}{2}}\!\!\!\left(\log{\frac{r_1}{r_0}}\right)^{\frac{1}{2}}+\frac{32}{3}t^3+3E_h^{1/2}\\
&\le&\!\!\!\! CE_h^{1/2}\left(\log{\frac{1}{s}}\right)^{1/2}+\frac{32}{3}t^3.
\label{eq:}
\end{eqnarray}
Plugging this into (\ref{eq:LB2CS}) yields the claim.\hfill$\square$

Next we establish a lower bound on the energy in the case that $w'$ remains for a long time in the well $W_r^+$ up to the boundary.
\begin{lem}\label{LB2:5}
Assume (\ref{eq:SetLB2}) and let $h\le\sqrt{\delta}$ and $w'(r)\in W_r^+$for all $r\in[\frac{\sqrt{\delta}}{8},1].$ Then the following bound holds
\begin{eqnarray}
\delta\ls\delta^{1/4}(\log{\frac{1}{h}})^{1/4}E_h^{1/4}+\frac{1}{\delta^{1/2}h^{1/4}}E_h^{1/2}+\frac{1+\log\frac{1}{\delta}}{\delta^{3/2}h}E_h.
\label{eq:}
\end{eqnarray}
\end{lem}
{\bf{Proof.}}\\
For every $l>0$ and $X\le 2l$ the elementary inequality
\begin{eqnarray}
X\le\frac{1}{16l^3}((16l^3X-4l^2X^2)+X^2(4l-X)^2)
\label{eq:}
\end{eqnarray}
holds. At first this looks fancy, but bringing the $X$ to the right hand side and multiplying by the factor $16l^3$ yields a polynomial in $X$ (for a fixed $l>0$): 
\begin{eqnarray}
P(X)=X^4+12l^2X^2-8lX^3.
\label{eq:}
\end{eqnarray}
By factorization we get
\begin{eqnarray}
P(X)=X^2(X-2l)(X-6l).
\label{eq:}
\end{eqnarray}
We have chosen the zeros of this polynomial in such away that $0\le P(X)$ for all $X\le 2l.$\\

If $r\in[\frac{\sqrt{\delta}}{8},1]$ then $w'(r)\in W_r^+$ and we can apply the above inequality with $l=r$ and $X=2r-w'(r)$ to obtain
\begin{eqnarray}
2r-w'\le\frac{1}{4r}(4r^2-w'^2)+\frac{1}{16r^3}(4r^2-w'^2)^2.
\label{eq:}
\end{eqnarray}
Similar to the proof of the lower bound for the e-cone we need a dyadic decomposition of the interval. Take $N\in\N$ s.t. $a=2^{-N}$ an $a$ satisfies $2^2a\le\sqrt{\delta}\le2^3a.$ Further let $a_j=2^{-j}.$ Then we start with the calculation

\begin{eqnarray}
\delta\!\!&=&\!\!\int\limits_0^1{1-w'\;dr}\\
&=&\!\!\int\limits_0^a{1-w'\;dr}+\int\limits_a^1{1-w'\;dr}\\
&\le&\!\! a+\int\limits_0^a{|w'|\;dr}+\int\limits_a^1{1-2r\;dr}+\int\limits_a^1{2r-w'\;dr}\\
&\le&\!\! a^2+\int\limits_0^a{|w'|\;dr}+\int\limits_a^1{2r-w'\;dr}\\
&\le&\!\!\frac{\delta}{16}+\int\limits_0^a{|w'|\;dr}+\frac{2}{\delta^{1/2}}\left|\int\limits_a^1{4r^2-w'^2\;dr}\right|+\frac{32}{\delta^{3/2}}\int\limits_a^1{|4r^2-w'^2|^2\;dr}\;\;\;\;\;\;\;\;\\
&=&\!\!\frac{\delta}{16}+(i)+|(ii)|+(iii).
\label{eq:}
\end{eqnarray}
Now we consider these three terms independently:
\begin{enumerate}
	\item To estimate the first one we use the inequalities from lemma \ref{LB2:4} with $s= \frac{h}{4}$ and $t=\frac{\sqrt{\delta}}{4}$
\begin{eqnarray}
\int\limits_0^a{|w'|\;dr}&\le&\int\limits_0^\frac{h}{4}{|w'|\;dr}+\int\limits_\frac{h}{4}^\frac{\sqrt{\delta}}{4}{|w'|\;dr}\\
&\le&\frac{\delta}{2\sqrt{3}}+\frac{1}{4}E_h^{1/2}+C\delta^{1/4}(\log{\frac{1}{h}})^{1/4}E_h^{1/4}.
\label{eq:}
\end{eqnarray}
		\item The second part we treat by dyadic decomposition. For this we calculate
\begin{eqnarray}
(ii)&=&\frac{2}{\delta^{1/2}}\int\limits_a^1{4r^2-w'^2\;dr}=\frac{2}{\delta^{1/2}}\sum\limits_{j=1}^N\int\limits_{a_j}^{a_{j-1}}{4r^2-w'^2\;dr}\\
&=&\frac{2}{\delta^{1/2}}\sum\limits_{j=1}^N{g_{a_j}(a_{j-1})-g_{a_j}(a_j)}.
\label{eq:LB2:(ii)}
\end{eqnarray}
By (\ref{eq:LB2:3.2}) we obtain
\begin{eqnarray}
|g_{a_j}(a_{j-1})-g_{a_j}(a_j)|\le\osc{I_{a_j}}{g_{a_j}}&\ls&\left(1+\frac{a_j^{1/2}}{h^{1/4}}\right)E_h^{1/2} \label{eq:LB2:(ii).1}\\
\ls\left(\frac{a_j}{h}\right)^{1/4}E_h^{1/2},
\label{eq:LB2:(ii).2}
\end{eqnarray}
where we have used the simple estimates $a_j\ge a\ge\frac{\sqrt{\delta}}{8}\ge\frac{h}{8}.$ Now taking the absolute value in (\ref{eq:LB2:(ii)}), using the triangle inequality and combining it with the estimate (\ref{eq:LB2:(ii).1}-\ref{eq:LB2:(ii).2}) yields
\begin{equation}
|(ii)|\ls\frac{1}{\delta^{1/2}h^{1/4}}E_h^{1/2}\sum\limits_{j=1}^N a_j^{1/4}\ls\frac{1}{\delta^{1/2}h^{1/4}}E_h^{1/2},
\label{eq:}
\end{equation}
where the sum is finite, since $\sum_{j=1}^N a_j^{1/4}\le\sum_{j=1}^\infty 2^{-j/4}<\infty.$

	  \item The third part can be treated similar to the one before. We start again with the decomposition 
\begin{eqnarray}
(iii)=\frac{32}{\delta^{3/2}}\int\limits_a^1{|4r^2-w'^2|^2\;dr}=\frac{32}{\delta^{3/2}}\sum\limits_{j=1}^N\int\limits_{I_{a_j}}{|4r^2-w'^2|^2\;dr}.
\label{eq:}
\end{eqnarray}
Now we apply the estimate (\ref{eq:LB2:3.1}), again from lemma \ref{LB2:3} to each summand. Hence
\begin{eqnarray}
(iii)\le\frac{32}{\delta^{3/2}}\sum\limits_{j=1}^N\left(\frac{1}{a_j}+\frac{a_j}{h}\right)E_h
\ls\frac{1+N}{\delta^{3/2}h}E_h,
\label{eq:}
\end{eqnarray}
where we have used $a_j\ge a\ge\frac{\sqrt{\delta}}{8}\ge\frac{h}{8}.$ The same chain of inequalities now implies $N\ls 1+\log\frac{1}{\delta}$ and the final bound
\begin{eqnarray}
(iii)\ls\frac{1+\log\frac{1}{\delta}}{\delta^{3/2}h}E_h.
\label{eq:}
\end{eqnarray}
\end{enumerate}
Together this yields
\begin{eqnarray}
\delta\ls\delta^{1/4}(\log{\frac{1}{h}})^{1/4}E_h^{1/4}+\frac{1}{\delta^{1/2}h^{1/4}}E_h^{1/2}+\frac{1+\log\frac{1}{\delta}}{\delta^{3/2}h}E_h.
\label{eq:}
\end{eqnarray}
\hfill$\square$\\
Similar we show a bound in the case that $w'$ stays in $W_r^-$ during the bulk. Such configurations should posses high energies.
\begin{lem}\label{LB2:6}
Assume (\ref{eq:SetLB2}) and let $h\le\sqrt{\delta}$ and $w'(r)\in W_r^-$ for all $r\in[\frac{\sqrt{\delta}}{8},1].$ Then the following bound holds
\begin{eqnarray}
\min\left\{1,\frac{1}{\delta\log\frac{1}{h}}\right\}\ls E_h.
\label{eq:}
\end{eqnarray}
\end{lem}
{\bf{Proof.}}\\ 
Again we start with the simple equation
\begin{eqnarray}
1-\delta=\int\limits_0^1{w'\;dt}=\int\limits_0^{\frac{\sqrt{\delta}}{8}}{w'\;dt}+\int\limits_{\frac{\sqrt{\delta}}{8}}^1{w'\;dt}.
\label{eq:BULK211}
\end{eqnarray}
It will be enough to estimate both integrals from above. For the second note that we have $w'(r)\le-\frac{3}{2}r$ for all $r\in[\frac{\sqrt{\delta}}{8},1]$
\begin{equation}
\int\limits_{\frac{\sqrt{\delta}}{8}}^1{w'\;dt}\le-\frac{3}{4}(1-\frac{\delta}{64}).
\label{eq:}
\end{equation}
The first one we can treat as $(i)$ in lemma \ref{LB2:5}. This time we use the inequalities from lemma (\ref{LB2:4}) with $s= \frac{h}{8},$ $t=\frac{\sqrt{\delta}}{8}$
\begin{eqnarray}
\int\limits_0^{\frac{\sqrt{\delta}}{8}}{|w'|\;dt}
&\le&\int\limits_0^\frac{h}{8}{|w'|\;dr}+\int\limits_\frac{h}{8}^\frac{\sqrt{\delta}}{8}{|w'|\;dr}\\
&\le&\frac{\delta}{8\sqrt{3}}+\frac{1}{8}E_h^{1/2}+C\delta^{1/4}(\log{\frac{1}{h}})^{1/4}E_h^{1/4}.
\label{eq:}
\end{eqnarray}
Starting from (\ref{eq:BULK211}) and using the estimates for the two terms above yields
\begin{equation}
1-\delta\le-\frac{3}{4}+\delta\left(\frac{3}{256}+\frac{1}{8\sqrt{3}}\right)+\frac{1}{8}E_h^{1/2}+C\delta^{1/4}(\log{\frac{1}{h}})^{1/4}E_h^{1/4}.
\label{eq:}
\end{equation}
Rearranging and $\delta\le1$ gives
\begin{equation}
1\ls E_h^{1/2}+\delta^{1/4}(\log{\frac{1}{h}})^{1/4}E_h^{1/4}.
\label{eq:}
\end{equation}
which proves the statement.\hfill$\square$\\

Now we can put things together to get the desired bound.
\begin{thm}\label{LB2}
Assume (\ref{eq:SetLB2}). Then the following lower bound on the energy holds
\begin{eqnarray}
h^2+\delta^{3/2}h^{3/2} \ls E_h.
\label{eq:}
\end{eqnarray}
\end{thm}
{\bf{Proof.}}\\
We start with the case $\delta\le h^{1/3}:$ Then $\delta^{3/2}h^{3/2}\le h^2$ and lemma \ref{LB2:SpB} implies
\begin{equation}
h^2+\delta^{3/2}h^{3/2}\ls h^2 \ls E_h.
\label{eq:}
\end{equation}
Now assume $\delta> h^{1/3}.$ 
Again consider the two different cases $\tau\ge\frac{\sqrt{\delta}}{8}$ and $\tau<\frac{\sqrt{\delta}}{8}.$ In the first case we can apply lemma (\ref{LB2:2}) like above and obtain
\begin{equation}
E_h\gs\min\{\tau^6,\tau^{3}h^{3/2}\}=\min\{\delta^{3},\delta^{3/2}h^{3/2}\}= \delta^{3/2}h^{3/2}.
\label{eq:}
\end{equation}
In the second case by definition of $\ta$ we know that the unique continuous representative $w'(r)\in W_r$ during the bulk $r\in[\frac{\sqrt{\delta}}{8},1].$\footnote{For a short proof of the existence and uniqueness of this continuous representative and further information about the set $\mathcal{A}_\delta,$ see \autocite[Lemma 19]{COT15}.}
Moreover we know by continuity that $w'(r)\in W_r^+$ or $w'(r)\in W_r^-$ during the bulk. If $w'(r)\in W_r^+$ for all $r\in[\frac{\sqrt{\delta}}{8},1]$ we apply Lemma (\ref{LB2:5}) and get
\begin{eqnarray}
E_h&\gs&\min\left\{\delta^3h^{1/2},\frac{\delta^3}{\log\frac{1}{h}},\frac{\delta^{5/2}h}{1+\log\frac{1}{\delta}}\right\}\\
&\gs&\min\left\{\delta^3h^{1/2},\frac{\delta^{5/2}h}{1+\log\frac{1}{\delta}}\right\}.
\label{eq:}
\end{eqnarray}
Else, $w'(r)\in W_r^-$ for all $r\in[\frac{\sqrt{\delta}}{8},1]$ we use lemma (\ref{LB2:6})
\begin{eqnarray}
 E_h\gs\min\left\{1,\frac{1}{\delta\log\frac{1}{h}}\right\}.
\label{eq:}
\end{eqnarray}
Comparing all of this cases leads to 
\begin{eqnarray}
 E_h\gs\min\left\{1,\frac{1}{\delta\log\frac{1}{h}},\delta^3h^{1/2},\frac{\delta^{5/2}h}{1+\log\frac{1}{\delta}},\delta^{3/2}h^{3/2}\right\}.
\label{eq:}
\end{eqnarray}
By $\delta>h^{1/3}$ it follows $\delta^3h^{1/2}\ge\delta^{3/2}h^{3/2}$ and 
\begin{eqnarray}
 E_h\gs\delta^{3/2}h^{3/2}
\label{eq:}
\end{eqnarray}
which proves the theorem.\hfill$\square$